\documentclass[12pt,a4paper,english,nofootinbib]{revtex4-2}
\usepackage{lmodern}
\usepackage{lmodern}

\usepackage[T1]{fontenc}
\usepackage[latin9]{inputenc}
\usepackage{babel}
\usepackage{calc}
\usepackage{amsmath}
\usepackage{amssymb}
\usepackage{graphicx}
\usepackage{esint}
\usepackage[unicode=true,pdfusetitle,
 bookmarks=true,bookmarksnumbered=false,bookmarksopen=false,
 breaklinks=false,pdfborder={0 0 1},backref=false,colorlinks=false]
 {hyperref}
\usepackage{xcolor}
\makeatletter

\@ifundefined{textcolor}{}
{%
 \definecolor{BLACK}{gray}{0}
 \definecolor{WHITE}{gray}{1}
 \definecolor{RED}{rgb}{1,0,0}
 \definecolor{GREEN}{rgb}{0,1,0}
 \definecolor{BLUE}{rgb}{0,0,1}
 \definecolor{CYAN}{cmyk}{1,0,0,0}
 \definecolor{MAGENTA}{cmyk}{0,1,0,0}
 \definecolor{YELLOW}{cmyk}{0,0,1,0}
}

\usepackage{babel}

\usepackage{graphicx}
\def\b{\begin{equation}}
\def\e{\end{equation}}

\@ifundefined{textcolor}{}{%
 \definecolor{BLACK}{gray}{0}
 \definecolor{WHITE}{gray}{1}
 \definecolor{RED}{rgb}{1,0,0}
 \definecolor{GREEN}{rgb}{0,1,0}
 \definecolor{BLUE}{rgb}{0,0,1}
 \definecolor{CYAN}{cmyk}{1,0,0,0}
 \definecolor{MAGENTA}{cmyk}{0,1,0,0}
 \definecolor{YELLOW}{cmyk}{0,0,1,0}
 }

\usepackage{latexsym}\usepackage{bm}

\begin{document}

\title{ A New Geometric Flow on 3-Manifolds: the $K$-Flow}
\author{{\normalsize{}{}{}{}{}{}{}{}{}{}{}{}{}{}{}{}{}Kezban Tasseten}}
\email{kezbantasseten@gmail.com}

\author{{\normalsize{}{}{}{}{}{}{}{}{}{}{}{}{}{}{}{}{}Bayram
Tekin}}
\email{btekin@metu.edu.tr}

\affiliation{Department of Physics, Middle East Technical University, 06800, Ankara,
Turkey}
\date{{\normalsize{}{}{}{}{}{}{}{}{}{}{{}{}{}{}\today}}}

\begin{abstract}
\noindent We define a new geometric flow, which we shall call the $K$-flow, on 3-dimensional Riemannian manifolds; and study the behavior of Thurston's model geometries under this flow both analytically and numerically. As an example, we show that an initially arbitrarily deformed homogeneous 3-sphere flows into a round 3-sphere and shrinks to a point in the unnormalized flow; or stays as a round 3-sphere in the volume normalized flow. The $K$-flow equation arises as the gradient flow of a specific purely quadratic action functional that has appeared as the quadratic part of New Massive Gravity in physics; and a decade earlier in the mathematics literature, as a new variational characterization of three-dimensional space forms. We show the short-time existence of the $K$-flow using a DeTurck-type argument.  
\end{abstract}
\maketitle
\[
\]
\tableofcontents{} 
\[
\]
\section{Introduction}

The Ricci flow equation, and its various modifications, arouse in the physics and mathematics literature almost simultaneously for different reasons. Let us discuss these briefly before we introduce our new flow. Friedan \cite{Friedan1,Friedan2} studied the extensions of a non-linear, 2-dimensional sigma model from a homogenous space to a compact Riemannian manifold $\mathcal{M}$. In the non-linear sigma model setting, the Riemannian metric $g_{i j}$ appears as a coupling ``constant''; and in quantum field theory, like every other parameter in the theory, it ``runs'' with the renormalization scale $\Lambda$. The continuum limit is achieved for $\Lambda \rightarrow \infty$, so $\Lambda^{-1}$ should be considered as a short-distance cut-off. Then, there is the so-called $\beta$ function, which in this case, is a tensor to be computed in perturbation theory.  To understand the results of the first and second-order perturbation theory, let us write the action of the non-linear sigma model in $2+\epsilon$ dimensions:
\begin{equation}
S[\phi^i] :=\frac{ \Lambda^\epsilon}{2}  \int d^{2+\epsilon}x\,  T^{-1}g_{i j}(\phi^k)\partial_\mu \phi^i \partial^\mu \phi^j, \label{nonlinear} 
\end{equation}
where $T$ is the temperature, and  $T^{-1}g_{i j}(\phi^k)$ is the positive-definite metric on (the target manifold) $\mathcal{M}$ which is assumed to be of a finite-dimensional Riemannian manifold. In quantum field theory, the running of the metric coupling, as one changes the renormalization scale, appears as a vector field in the infinite-dimensional space of Riemannian metrics on $\mathcal{M}$; and for the model (\ref{nonlinear}), at low-temperature expansion (which corresponds to loop expansion in perturbative quantum field theory),  one has \cite{Friedan1,Friedan2}
\begin{equation}
\beta_{ij} (\Lambda):=-\Lambda^{-1}\frac{\partial}{ \partial \Lambda^{-1}} \left(T^{-1}g_{ij}(\Lambda)\right)=-\epsilon T^{-1}g_{ij} +R_{i j} +\frac{T}{2} R_{i pqr} R_j\,^{pqr} + \mathcal{O}(T^2).
\end{equation}
In the $\epsilon \rightarrow 0$ limit, at the lowest order, and defining $t:= \log \Lambda$, and $T= -2$, one arrives at the Ricci-flow equation
\begin{equation}
\partial_t g_{i j} = -2 R_{i j}.
\end{equation}
Independently, this equation was introduced by Hamilton \cite{Hamilton1} to deform metrics along their Ricci tensor on a given manifold, especially a 3-manifold to prove the Poincar\'e conjecture about $S^3$: that is every simply connected, closed 3-manifold is homeomorphic to $S^3$. Hamilton achieved a great deal with the Ricci flow technique: he proved the short-time existence of the flow for any initial metric and classified closed 3-manifolds with positive Ricci curvature. But the final proof of the Poincar\'e conjecture and the more general Thurston's uniformization conjecture \cite{Thurston} about 3 manifolds was given, rather famously, by Perelman \cite{Perelman1,Perelman2} who used a non-minimally coupled scalar-tensor action to derive a modified version of Ricci-flow and surgery techniques.  For a detailed exposition of these papers, see \cite{Lott}.

From the physics point of view, Ricci flow or some geometric flow of the type 
 \begin{equation}
\partial_t g_{i j} = \Phi_{i j},
\end{equation}
where  $\Phi_{i j}$ is some tensor built on the metric tensor and its derivatives is a great tool to understand the transitions between the saddle points (i.e those metrics satisfying 
$ \Phi_{i j}=0$) of a gravity theory. If the flow equation comes from the Euclidean action of the theory, these transitions will most probably correspond to some quantum mechanical effects in the Lorentzian version of the theory. Therefore, understanding these flows is quite important: for this reason in \cite{Cotton1} using the third order Cotton tensor, which is an obstruction to conformal flatness in 3 dimensions, was used to define the Cotton flow and the model geometries under this flow were studied. The discussion was extended in \cite{Cotton2} where a refinement on the definition of Cotton solitons was given.

The $K$-flow that we shall introduce in this work is based on a physical model of gravity that emerged in \cite{nmg1,nmg2} and has been extensively studied in the literature as the first example of a non-linearly massive gravity theory \cite{Gullu1,Gullu2} albeit in $2+1$ dimensions.  The action of the model is given (in the Lorentzian signature) as
\begin{equation}
S[g]:=-\frac{1}{\kappa^{2}}\int d^{3}x\,\sqrt{-g}\left(R-2\Lambda_{0}+\frac{1}{M^{2}}K\right), \hskip 1 cm K :=R_{\mu\nu}R^{\mu\nu}-\frac{3}{8}R^{2},
\end{equation}
of which the source-free field equations are 
\begin{equation}
G_{\mu\nu}+\Lambda_{0}g_{\mu\nu}-\frac{1}{2M^{2}}K_{\mu\nu}=0.\label{eq:NMG_eom}
\end{equation}
$\Lambda_0$ and $M$ are the bare cosmological constant and the mass parameter; and the $K$-tensor reads
\begin{equation}
K_{\mu\nu}=2\square R_{\mu\nu}-\frac{1}{2}\left(\nabla_{\mu}\nabla_{\mu}+g_{\mu\nu}\square\right)R+4R_{\mu\alpha\nu\beta}R^{\alpha\beta}-\frac{3}{2}RR_{\mu\nu}-g_{\mu\nu}K,
\label{K-tensor}
\end{equation}
with the trace $K=g^{\mu\nu}K_{\mu\nu}$, and $\square := g^{\mu \nu}\nabla_\mu \nabla_\nu$. 

Unlike the 3-dimensional Einstein's gravity based on the Einstein-Hilbert action, the theory is locally non-trivial: it has a massive spin-2 graviton (and no other degrees of freedom) with the mass around its (anti)-de Sitter background given as 
\begin{equation}
m_{g}^{2}=M^2-\frac{\Lambda}{2}, \hskip 2 cm \Lambda = 2 M^2 \Big ( -1 \pm \sqrt{1 - \frac{\Lambda_0}{M^2}} \,\, \Big ). 
\end{equation}
Of course, in addition to the perturbative sector, the theory also has many non-perturbative solutions, such as black holes, solitons, etc. So it is a rather rich theory of gravity that deservedly attracted a lot of attention in physics.
If one drops the Ricci scalar and the bare cosmological constant in the action, one arrives at the purely quadratic theory with a single massless degree of freedom \cite{Deser, Gullu2}. 
\begin{equation}
S_K:=-\frac{1}{M^2\kappa^{2}}\int d^{3}x\,\sqrt{-g}\left(R_{\mu\nu}R^{\mu\nu}-\frac{3}{8}R^{2}\right).
\label{K-theory}
\end{equation}
As it sometimes happens, this action, before it appeared as a model of 3D gravity, was studied in the mathematics literature \cite{Viac} from a different vantage point which we shall discuss briefly below. From the physical point of view, at short distances (at the so-called ultra-violet regime), higher curvature theories are much better behaved. Therefore, in what follows, we shall build a geometric flow which will be a gradient flow of the Riemannian version (that is with a positive-definite metric) of the action (\ref{K-theory}).

Before we move on, let us state our conventions: we will have to work in both local coordinates, orthonormal frame, and in the index-free form of the tensors. Therefore we will not shy away from calling the index-free Riemann tensor as $\text{Riem}$ and the Ricci tensor as $\text{Ric}$. In the coordinate-adapted basis, we shall use the indices ${i, j, k,l...}$, while in the orthonormal basis, we shall use the indices ${a,b,c,...}$. We will always work with positive-definite Riemannian metric and our sign convention for the Riemann tensor follows from  $[\nabla_i, \nabla_j] w_k = R_{i j k}\,^l w_l$, and the Ricci tensor is $R_{i j}:= R^{k}\,_{i k j}$. 

\subsection{The $K$-flow}

Let ${\mathcal{M}}$ be a 3-dimensional Riemannian manifold with a (positive definite) metric ${\bf g}$ on it; and let us introduce a deformation parameter $t$, not a coordinate on the manifold, akin to the renormalization group parameter introduced above; and let {\bf K} be the symmetric $K$ tensor whose components in local coordinates are given as the negative of (\ref{K-tensor}). Then the $K$-flow is defined as\footnote{For all intents and purposes, we can envisage a one-parameter family of Riemannian manifolds $( {\mathcal{M}}(t),  {\bf g}(t) )$, and the flow smoothly deforms one such member to an infinitesimally close one as the ``time'' $t$ evolves. With this vantage point, one can see $ {\mathcal{M}}(t)$  as a spacelike hypersurface in a spacetime with $t$ being the time dimension, and the evolution is the flow equation instead of the usual hyperbolic equation coming from Einstein's equation.  }
		\begin{eqnarray}  
		\partial_{t}{\bf g}(t)&=&\alpha {\bf K}(t) \qquad {\bf g}(t=0)={\bf g}_{0}, 
	\end{eqnarray}
where $\alpha$ is a constant later to be specified.  The flow of the inverse metric easily follows as $\partial_{t}{\bf g}^{-1}(t)=-\alpha {\bf K}^{-1}(t)$.
The $K$-tensor is slightly cumbersome, it is better to split it into two parts as it was introduced in \cite{Tekin_bulk}.  For this purpose let us introduce some local coordinates and define two relevant tensors:
\begin{equation}
		C_{ij}:=\eta^{\phantom{i}kl }_{i}\nabla_{k}S_{lj} \qquad S_{ij}:=R_{ij}-\frac{1}{4}g_{ij}R,
	\end{equation}
where the $C_{i j}$ is the Cotton tensor and $S_{i j}$ is the Schouten tensor, while $\eta^{i kl }$ is the completely anti-symmetric tensor. Using these, we define a traceless, 3D Bach-like tensor as  
\begin{eqnarray}
		H_{ij}&:=& \frac{1}{2}\eta^{\phantom{i}kl }_{i}\nabla_{k}C_{lj}+\frac{1}{2}\eta^{\phantom{j}kl }_{j}\nabla_{k}C_{li} 
		=3S_{\phantom{k}i }^{k}S_{kj}-g_{ij}S_{kl}S^{kl}-\Delta S_{ij}+\nabla_{i}\nabla_{j}S, \label{Bach}
	\end{eqnarray}
where  $ \Delta=g^{ij}\nabla_{i}\nabla_{j}=\nabla^{i}\nabla_{i}$, $S= R/4$ and the trace part of $K$ as 
\begin{eqnarray}
		J_{ij}&:=& \frac{1}{2}\eta^{\phantom{i}kl }_{i}\eta^{\phantom{j}mn }_{j}S_{km}S_{ln} =S_{\phantom{k}i }^{k}S_{kj}-SS_{ij}+\frac{1}{2}g_{ij}(S^{2}-S_{kl}S^{kl}).
\label{Jtensor}	
	\end{eqnarray}
	Therefore we have 
	\begin{equation} 
		K_{i j}=2J_{i j}+2H_{i j}, \hskip 1 cm  K=2J =   \dfrac{3}{8} R^{2}- R_{i j} R^{ ij}.  
		\label{ikinciK}
	\end{equation}
	Important remark: Let us repeat that this $K_{i j}$ is negative of the $K_{\mu \nu}$ given in  (\ref{K-tensor}). From now on, we will work with (\ref{ikinciK}).
	
We will study the flows of model geometries in the orthonormal frame. To define our tensor $K$ in the orthonormal frame, we first define the $J$ and $H$ one-forms as
	\begin{equation}
		J^{a}=\frac{1}{2}\epsilon_{\phantom{a}bc}^{a}\star (S^{b}\wedge S^{c}) \qquad H^{a}=\star D\star C^{a},
	\end{equation}
	where $S^{a}$ is the Schouten one-form and $\star C^{a}$ is the Cotton one-form; and the star denotes the Hodge dual as usual.  More explicitly, we have the following: Let $w^a\,_b$ be the connection one-forms, and $D:= d  +\omega \wedge$ be the covariant derivative one-form such that $D e^a=0$. Then the curvature 2-form reads $R^a\,_b:= d \omega^a\,_b + \omega^a\,_c \wedge \omega^c\, _b= D\omega^a\,_b$, from which follows the Rici-one form $Ric_a :=\iota_b R^b\,_a$ by interior product; and the scalar curvature $R :=\iota_a Ric^a$. The Schouten one-form (in 3D) then is defined as $S^{a}:= Ric^a -\frac{1}{4} R e^a$, from which follows the Cotton 2-form $C^a= DS^a$.
	Then we have the $K$ one-form:
	\begin{equation}
	K^a := \epsilon_{\phantom{a}bc}^{a}\star (S^{b}\wedge S^{c}) + 2\star D\star C^{a},
\end{equation}
and the $K$ flow equation is 
\begin{equation}
\partial_t e^a = \alpha K^a.
\end{equation}
In the table, we give a list of the homogenous geometries and their Ricci and $K_{ij}$ tensors.\footnote{See \cite{Cotton2} for a similar table where the Cotton tensors of the Thurston geometries were given.} From this table, it should be clear to the reader that working in these coordinate forms, generically, would not yield autonomous differential equations for most of the geometries. Therefore, following similar works in the Ricci flow \cite{Isenberg,KM,CK}, and the Cotton flow {\cite{Cotton1}, we will work in the Milnor frame \cite{Mil} in which one obtains autonomous nonlinear ODEs.

\begin{center}
\resizebox{\columnwidth}{!}{%
    \begin{tabular}{ | l | l | l | l| l | p{5 cm}|}
    \multicolumn{4}{c}{A list of Thurston geometries and their Ricci and the $K$ curvatures} \\
\cline{1-4}
    Geometry &\qquad \qquad Metric & $\qquad R_{ij}$ & $\qquad K_{ij}$  \\ \hline
    $\mathbb{R}^3$ & $ds^2=dx^2+dy^2+dz^2$&$ R_{ij}=0$& $K_{ij}=0$\\ \hline
   $S^3$ & $ds^2=dx^2+\sin^2 x dy^2+ (dz+ \cos x\, dy)^2$& $R_{ij}=\frac{1}{2}g_{ij}$ & $K_{ij}=\frac{1}{32} g_{i j}$\\ \hline
   $\mathbb{H}^3$ & $ds^2=\frac{1}{x^2} (dx^2+dy^2+dz^2)$& $ R_{ij}=-2 g_{ij}  $ & $K_{ij}=\frac{1}{2}g_{ij}$ \\ \hline
   $S^2 \times \mathbb{R}$ & $ds^2=dx^2+\sin^2 x \, dy^2+dz^2$ &$ R_{11}=1, \, R_{22}=\sin^2 x $ &  $K_{11}=-\frac{1}{2}, K_{22}= -\frac{\sin^2x}{2},K_{33}=\frac{1}{2}$ \\ \hline
   $\mathbb{H}^2 \times \mathbb{R}$ & $ds^2=\frac{1}{x^2} (dx^2+dy^2)+dz^2$ &$ R_{11}=R_{22}=-\frac{1}{x^2} $  &  $K_{11}= K_{22}=-\frac{1}{2 x^2}, K_{33}=\frac{1}{2}$ \\ \hline
   $Sol$ & $ds^2=e^{2 z}dx^2+e^{-2 z}dy^2+dz^2$ & $R_{33}=-2$ & $ K_{11}= K_{22}=\frac{5}{2}e^{2 z}, K_{33}=-\frac{15}{2} $  \\ \hline
   $Nil$ & $ds^2=dx^2+dy^2+(dz+x \,dy )^2$ & $R_{11}=-\frac{1}{2}, \, R_{22}=\frac{x^2-1}{2}, $ & $K_{11}=-\frac{63}{32}, K_{22}=\frac{21}{32}\left(5 x^2-3\right)$, \\
   & & $R_{33}=\frac{1}{2},\, R_{23}=\frac{x}{2}$ & $K_{33}=\frac{105}{32}, \, K_{23}=\frac{105x}{32}$ \\ \hline
   $SL(2, \mathbb{R})$ & $ds^2=\frac{1}{x^2}  (dx^2+dy^2)+(dz+\frac{1}{x} d y )^2$ & $ R_{11}=-\frac{3}{2x^2},\, R_{22}=-\frac{1}{x^2},$  & $K_{11}=-\frac{159}{32 x^2},\, K_{22}=\frac{41}{16 x^2}, $ \\ 
   & & $ R_{23}=\frac{1}{2x},R_{33}=\frac{1}{2}$  & $ K_{23}=\frac{241}{32 x},\, K_{33}=\frac{241}{32}$\\
   \hline
    \end{tabular}%
    }
\end{center}
\vskip 0.7 cm 
We remark that although there are eight geometries in the table as in Thurston's geometrization conjecture, from classification of Milnor one more unimodular simply-connected Lie group, $\widetilde{Isom}(R^{2})$ appears which collapses to the same maximally symmetric geometry  of  $\mathbb{R}^3$. 

\section{The short-time Existence of the $K$-Flow}

Before we move on to the study of model geometries we would like to discuss the existence of the flow. The procedure we will follow was suggested by  DeTurck \cite{DeTurck}  in giving concise proof of the short-time existence of the Ricci flow. As an evolutionary partial differential equation, Ricci flow is not strictly parabolic. DeTurck used the method of principal symbols of linear differential operators to prove the short-time existence of the Ricci flow. We will follow the same procedure which amounts to evaluating the tensor for small perturbations.  All of the high-frequency, short-wavelength modes are expected to diffuse for a parabolic equation. For that purpose, the linearized version of the $K$-tensor around a flat background is sufficient. But, first, for the sake of generality, let us find the linearized version of this tensor about a generic background geometry that is not necessarily flat.

\subsection{Variation of the $K$-Tensor about a Background}
We vary our metric $g$ in around a background as $g=\bar{g}+h$, where $\bar{g}$ is the background and $h$ is the variation, so that 
\begin{equation}
	    g_{i j}= \bar{g}_{i j}+h_{i j}, \qquad
		g^{i j}= \bar{g}^{i j}-h^{i j} + {\mathcal{O}}(h^2).
\end{equation}	 
We raise and lower indices by $\bar{g}$ and take the variation up to the second order.
$\delta(g_{ij})=h_{ij}$ and $\delta(g^{ij})=-h^{ij}$. Then we have the following variations for the terms that do not involve derivatives of the curvatures

\begin{eqnarray}
	\delta(\Gamma^{i}_{jk})&=&\frac{1}{2}\bar{g}^{il}(-\bar{\nabla}_{l}h_{jk}+\bar{\nabla}_{j}h_{kl}+\bar{\nabla}_{k}h_{jl}), \nonumber \\
		\delta(R_{ij})&=&\frac{1}{2}(-\bar{\nabla^{2} }h_{ij}+\bar{\nabla}_{k}\bar{\nabla}_{j}h_{\phantom{k}i }^{k}+\bar{\nabla}_{k}\bar{\nabla}_{i}h_{\phantom{k}j }^{k}-\bar{\nabla}_{i}\bar{\nabla}_{j}h), \nonumber \\
		\delta(R)&=& -h^{kl}\bar{R}_{kl}+\bar{\nabla}_{k}\bar{\nabla}_{l}h^{kl}-\bar{\nabla^{2} }h, \nonumber \\
		\delta(S_{\phantom{k}i }^{k}S_{kj})&=&-h^{kl}\bar{R}_{ki}\bar{R}_{lj}+\delta(R_{ki})\bar{R}_{\phantom{k}j }^{k}+\delta(R_{kj})\bar{R}_{\phantom{k}i }^{k}-\frac{1}{2}\delta(R_{ij})\bar{R}-\frac{1}{2}\delta(R)\bar{R}_{ij},\nonumber \\
		&+&\frac{1}{8}\delta(R)\bar{g}_{ij}\bar{R}+\frac{1}{16}h_{ij}\bar{R}^{2}, \nonumber \\
		\delta(SS_{ij})&=& \frac{1}{4}\delta(R_{ij})\bar{R}+\frac{1}{4}\delta(R)\bar{R}_{ij}-\frac{1}{8}\delta(R)\bar{g}_{ij}\bar{R}-\frac{1}{16}h_{ij}\bar{R}^{2}, \nonumber \\
		\delta \left(g_{ij}(S_{kl})^{2}\right)&=& 2\delta(R_{kl})\bar{g}_{ij}\bar{R}^{kl}-\frac{1}{2}\delta(R_{kl})\bar{g}^{kl}\bar{g}_{ij}\bar{R}-\frac{1}{8}\delta(R)\bar{g}_{ij}\bar{R},\nonumber \\
		&+&\frac{1}{2}h_{kl}\bar{g}_{ij}\bar{R}^{kl}\bar{R}-2h^{kn}\bar{g}_{ij}\bar{R}_{kl}\bar{R}_{\phantom{l}n }^{l}+h_{ij}(\bar{R}_{kl})^{2}-\frac{5}{16}h_{ij}\bar{R}^{2} \nonumber \\
		\delta(g_{ij}S^{2})&=&\frac{1}{8}\delta(R)\bar{g}_{ij}\bar{R}+\frac{1}{16}h_{ij}\bar{R}^{2}. 
\end{eqnarray} 		
The terms that involve the derivatives of the curvatures  vary as	
\begin{eqnarray}	
\delta(\nabla_{i}\nabla_{j}S)&=&\bar{\nabla}_{i}\bar{\nabla}_{j}\delta(S)-\frac{1}{8}(\bar{\nabla}_{i}h^{\phantom{j}k}_{j}+\bar{\nabla}_{j}h^{\phantom{i}k}_{i}-\bar{\nabla}^{k}h_{ij})\bar{\nabla}_{k}\bar{R},\nonumber \\
		\delta(\nabla^{2}  S_{ij})&=& \bar{\nabla^{2} }\delta(S_{ij})-h^{kl}\bar{\nabla}_{k}\bar{\nabla}_{l}\bar{S}_{ij}+\frac{1}{2}(\bar{\nabla}^{k}\bar{\nabla}^{l}h_{ki}-\bar{\nabla}_{k}\bar{\nabla}_{i}h^{kl}-\bar{\nabla^{2} }h_{\phantom{l}i}^{l})S_{lj} \nonumber \\
		&+&\frac{1}{2}(\bar{\nabla}^{k}\bar{\nabla}^{l}h_{kj}-\bar{\nabla}_{k}\bar{\nabla}_{j}h^{kl}-\bar{\nabla^{2} }h_{\phantom{l}j}^{l})S_{li} 
		+ \frac{1}{2}(\bar{\nabla}^{k}h-2\bar{\nabla}_{l}h^{kl})\bar{\nabla}_{k}\bar{S}_{ij}\nonumber \\ 
		&+&(\bar{\nabla}^{k}h_{\phantom{l}i}^{l}-\bar{\nabla}^{l}h_{\phantom{k}i}^{k}-\bar{\nabla}_{i}h^{kl})\bar{\nabla}_{l}\bar{S}_{kj}+(\bar{\nabla}^{k}h_{\phantom{l}j}^{l}-\bar{\nabla}^{l}h_{\phantom{k}j}^{k}-\bar{\nabla}_{j}h^{kl})\bar{\nabla}_{l}\bar{S}_{ki}.
\end{eqnarray} 
Finally, collecting all the pieces and suppressing the bar over the background tensors for notational simplicity, one finds the variation of the $ K$ tensor as 
\begin{eqnarray}
		\delta(K_{ij})&=& \Delta^{2}h_{ij}-\Delta \nabla_{k}\nabla_{i}h^{\phantom{j}k}_{j}-\Delta \nabla_{k}\nabla_{j}h^{\phantom{i}k}_{i}+\frac{1}{2}\nabla_{i}\nabla_{j}\nabla_{k}\nabla_{l}h^{kl}+\Delta \nabla_{i}\nabla_{j}h \nonumber \\
		&-&\frac{1}{2}\nabla_{i}\nabla_{j}\Delta h+\frac{1}{2} g_{ij}(\Delta\nabla_{k}\nabla_{l}h^{kl}-\Delta^{2}h) +g_{ij}(6h^{kl}R^{\phantom{k}n}_{k}R_{ln}-\frac{13}{4}h^{kl}R_{kl}R)\nonumber \\
		&+&g_{ij}\Bigl(3R_{kl}(\Delta h^{kl}-\nabla_{n}\nabla^{k}h^{ln}-\nabla_{n}\nabla^{l}h^{kn}+\nabla^{k}\nabla^{l}h)+\frac{13}{4}R(\nabla_{k}\nabla_{l}h^{kl}-\Delta h)\Bigr)\nonumber \\
		&+&g_{ij}(-\frac{1}{2}\Delta (h^{kl}R_{kl})-\frac{1}{2}\nabla_{k}R\nabla_{l}h^{kl}+\frac{1}{4}\nabla_{k}R\nabla^{k}h-\frac{1}{2}h^{kl}\nabla_{k}\nabla_{l}R) \nonumber \\
		&-&\frac{1}{2}\nabla_{i}\nabla_{j}(h^{kl}R_{kl})+\frac{9}{2}R_{ij}(\Delta h-\nabla_{k}\nabla_{l}h^{kl}) \nonumber \\
		&+&R_{ik}(3\nabla^{l}\nabla^{k}h_{lj}-3\Delta h^{\phantom{j}k}_{j}+5\nabla^{l}\nabla_{j}h^{\phantom{l}k}_{l}-4\nabla_{j}\nabla^{k}h)+i \leftrightarrow j \nonumber \\
		&+& \frac{9}{4}(\Delta h_{ij}-\nabla_{k}\nabla_{i}h^{\phantom{j}k}_{j}-\nabla_{k}\nabla_{j}h^{\phantom{i}k}_{i}+\nabla_{i}\nabla_{j}h)+2\nabla_{k}R_{ij}\nabla_{l}h^{kl}-\nabla_{k}R_{ij}\nabla^{k}h \nonumber \\
		&+&2\nabla_{k}R_{il}(\nabla^{k}h^{\phantom{j}l}_{j}-\nabla^{l}h^{\phantom{j}k}_{j}+\nabla_{j}h^{kl}) +i \leftrightarrow j \nonumber \\
		&+& \frac{1}{4}\nabla_{k}R(\nabla^{k}h_{ij}-\nabla_{i}h^{\phantom{j}k}_{j}-\nabla_{j}h^{\phantom{i}k}_{i})+2h^{kl}\nabla_{k}\nabla{l}R_{ij}+\frac{1}{2}h_{ij}\Delta R \nonumber \\
		 &-&3h_{ij}R_{kl}R^{kl}+\frac{9}{2}R_{ij}R_{kl}h^{kl}-8R_{ik}R_{jl}h^{kl}+\frac{13}{8}h_{ij}R^{2}.
\end{eqnarray}
The above computation was about a generic background space, let us consider the case for an Einstein (maximally symmetric in 3D) space for which $\bar{R}_{ij}=2 \Lambda \bar{g}_{ij}$
	\begin{eqnarray}		\delta(K_{ij})&=&\Delta^{2}h_{ij}-\Delta \nabla_{k}\nabla_{i}h^{\phantom{j}k}_{j}-\Delta \nabla_{k}\nabla_{j}h^{\phantom{i}k}_{i}+
		\frac{1}{2}\nabla_{i}\nabla_{j}\nabla_{k}\nabla_{l}h^{kl}+\Delta \nabla_{i}\nabla_{j}h \nonumber \\
		&-&\frac{1}{2}\nabla_{i}\nabla_{j}\Delta h+\frac{1}{2} g_{ij}(\Delta\nabla_{k}\nabla_{l}h^{kl}-\Delta^{2}h)+g_{ij}\Lambda \left(-\frac{3}{2}\nabla_{k}\nabla_{l}h^{kl}+\frac{1}{2}\Delta h+3\Lambda h\right) \nonumber \\
		&+&\Lambda\left(\frac{5}{2}\nabla^{k}\nabla_{i}h_{jk}+\frac{5}{2}\nabla^{k}\nabla_{j}h_{ik}+\frac{3}{2}\Delta h_{ij}-\frac{7}{2}\nabla_{i}\nabla_{j}h\right)-\frac{19}{2}\Lambda^{2}h_{ij}.
	\end{eqnarray}
Thus the variation of the $K$ tensor is dominated by the fourth-order derivatives of the perturbation. In the case of a flat background, $\Lambda=0$, we get the linearized version of the $K$ tensor ($K^L_{ij}$), namely 
\begin{eqnarray}
		K_{ij}^{L}&=&\Delta^{2}h_{ij}-\Delta \nabla_{k}\nabla_{i}h^{\phantom{j}k}_{j}-\Delta \nabla_{k}\nabla_{j}h^{\phantom{i}k}_{i}+\frac{1}{2}\nabla_{i}\nabla_{j}\nabla_{k}\nabla_{l}h^{kl}+\Delta \nabla_{i}\nabla_{j}h \nonumber \\
		&-&\frac{1}{2}\nabla_{i}\nabla_{j}\Delta h+\frac{1}{2} g_{ij}\left(\Delta\nabla_{k}\nabla_{l}h^{kl}-\Delta^{2}h\right).
\label{linearK}
\end{eqnarray}
In the jargon of physics, in the transverse (i.e. $\nabla_i h^{i j}=0$) and the traceless (i.e. $h=0$) gauge, in Cartesian coordinates, we just have $K_{ij}^{L}=\left (\partial_k \partial_k \right )^2 h_{ij}$ which is a bi-laplacian operator. This suggests that we have a fourth-order elliptic operator; but to show this properly, let us calculate the principal symbol.

\subsection{The principal symbol of the relevant operator}
Using (\ref{linearK}), we define 
\begin{eqnarray}
		[D(K_{g}(h))]_{ij}&:=&\frac{1}{2}\nabla_{i}\nabla_{j}\nabla_{k}\nabla_{l}h^{kl} +g^{mn}\Bigl(g^{pq}\nabla_{m}\nabla_{n}\nabla_{p}\nabla_{q}h_{ij}- \nabla_{m}\nabla_{n}\nabla_{k}\nabla_{i}h^{\phantom{j}k}_{j} \nonumber \\
		&-&\nabla_{m}\nabla_{n} \nabla_{k}\nabla_{j}h^{\phantom{i}k}_{i}+\nabla_{m}\nabla_{n}\nabla_{i}\nabla_{j}h-\frac{1}{2}\nabla_{i}\nabla_{j}\nabla_{m}\nabla_{n} h 
		 \nonumber \\
		&+&\frac{1}{2} g_{ij}\left( \nabla_{m}\nabla_{n}   \nabla_{k}\nabla_{l}h^{kl}-g^{pq}\nabla_{m}\nabla_{n} \nabla_{p}\nabla_{q}h \right)\Bigr).	
	\end{eqnarray}
Let $\zeta \in C^{\infty}(T^{*}{\mathcal{M}})$ be a covector. The {\it principal symbol} of the linear differential operator $D(K_{g})$ is the bundle homomorphism that is obtained by replacing the covariant derivative $\nabla_{k}$ by a covector $\zeta_{k}$. Hence we have 
	\begin{eqnarray}
		\Big[\hat{\sigma}[D(K_{g})](\zeta)(h)\Big]_{ij}&=&\frac{1}{2}\zeta_{i}\zeta_{j}\zeta_{k}\zeta_{l}h^{kl} +g^{mn}\Bigl( g^{pq}\zeta_{m}\zeta_{n}\zeta_{p}\zeta_{q}h_{ij}- \zeta_{m}\zeta_{n}\zeta_{k}\zeta_{i}h^{\phantom{j}k}_{j} \nonumber \\
		&-&\zeta_{m}\zeta_{n} \zeta_{k}\zeta_{j}h^{\phantom{i}k}_{i}+\zeta_{m}\zeta_{n}\zeta_{i}\zeta_{j}h-\frac{1}{2}\zeta_{i}\zeta_{j}\zeta_{m}\zeta_{n} h 
		+ \nonumber \\
		&+& \frac{1}{2} g_{ij}( \zeta_{m}\zeta_{n}   \zeta_{k}\zeta_{l}h^{kl}-g^{pq}\zeta_{m}\zeta_{n} \zeta_{p}\zeta_{q}h)\Bigr).
	\end{eqnarray}	
	Then for all covectors $\zeta\neq 0$ and for all $h_{ij}\neq 0$ we find the inner product:
	\begin{eqnarray}
		\langle\hat{\sigma}[D(K_{g})](\zeta)(h),h \rangle&=&\frac{1}{2}\zeta_{i}\zeta_{j}\zeta_{k}\zeta_{l}h^{kl} h^{ij} +g^{mn}\Bigl(g^{pq}\zeta_{m}\zeta_{n}\zeta_{p}\zeta_{q}h_{ij}- \zeta_{m}\zeta_{n}\zeta_{k}\zeta_{i}h^{\phantom{j}k}_{j}\nonumber \\
		&-&\zeta_{m}\zeta_{n} \zeta_{k}\zeta_{j}h^{\phantom{i}k}_{i}+\zeta_{m}\zeta_{n}\zeta_{i}\zeta_{j}h-\frac{1}{2}\zeta_{i}\zeta_{j}\zeta_{m}\zeta_{n} h  \nonumber \\
		&+&\frac{1}{2} g_{ij}( \zeta_{m}\zeta_{n}   \zeta_{k}\zeta_{l}h^{kl}-g^{pq}\zeta_{m}\zeta_{n} \zeta_{p}\zeta_{q}h)\Bigr) h^{ij}.
	\end{eqnarray}	
It is easy to see that for $h_{ij} =\zeta_{i}\zeta_{j}$, which is a pure gauge mode, the above equation vanishes, showing that the operator is not elliptic. To obtain an elliptic operator, devoid of these zero modes, we commute the covariant derivatives in the linearized $K$ tensor. Whenever we commute the covariant derivatives, we obtain the Riemann tensor, that is 
\begin{equation}
\Delta \nabla_{k}\nabla_{i}h^{\phantom{j}k}_{j}=\Delta \nabla_{i}\nabla_{k}h^{\phantom{j}k}_{j}+\Delta Riem \cdot h.
\end{equation}
Since  the only contribution to the principal symbol comes from the highest order derivatives, we can commute derivatives with ease to obtain 
\begin{eqnarray}
		[D(K_{g}(h))]_{ij}&=& \Delta^{2}h_{ij} + \frac{1}{2} g_{ij}(\Delta\nabla_{k}\nabla_{l}h^{kl}-\Delta^{2}h) \nonumber \\
		&+& \nabla_{i}\left (- \nabla_{k}\Delta h^{\phantom{j}k}_{j}+\frac{1}{4}\nabla_{j}\nabla_{k}\nabla_{l}h^{kl}+\frac{1}{4}\nabla_{j}\Delta h\right) \nonumber \\
		&+&\nabla_{j}\left (- \nabla_{k}\Delta h^{\phantom{i}k}_{i}+\frac{1}{4}\nabla_{i}\nabla_{k}\nabla_{l}h^{kl}+\frac{1}{4}\nabla_{i}\Delta h \right),
	\end{eqnarray}
which suggests that we define the following one-form 
\begin{equation}
V_{j} :=- \nabla_{k}\Delta h^{\phantom{j}k}_{j}+\frac{1}{4}\nabla_{j}\nabla_{k}\nabla_{l}h^{kl}+\frac{1}{4}\nabla_{j}\Delta h,
\end{equation}
 in terms of which  one has  
\begin{equation}
[D(K_{g}(h))]_{ij}= \Delta^{2}h_{ij} + \frac{1}{2} g_{ij}(\Delta\nabla_{k}\nabla_{l}h^{kl}-\Delta^{2}h)+\nabla_{i}V_{j}+\nabla_{j}V_{i}.
\end{equation}
We can express the one-form $V_{j}$ as 
	\begin{equation}
	V_{j}=\frac{1}{2}\nabla^{k}\left[D(Ric_{g})(h)\right]_{jk}-g_{jk}g^{mn}\Delta \left[D(\Gamma_{g})(h)\right]^{k}_{mn},
	\label{Veq}
	\end{equation}
where 
\begin{equation}
[D(Ric_{g})]:C^{\infty}(S_{2}T^{*}{\mathcal{M}})\rightarrow C^{\infty}(S_{2}T^{*}{\mathcal{M}})
\end{equation}
denotes the linearization of the Ricci tensor, $S_{2}$ refers to symmetric two tensors; and 
\begin{equation}
[D(\Gamma_{g})]:C^{\infty}(S_{2}T^{*}{\mathcal{M}})\rightarrow C^{\infty}(S_{2}T^{*}{\mathcal{M}}\times T{\mathcal{M}})
\end{equation}
denotes the linearization of the Levi-Civita connection. The connection is not a tensor, but we may define a vector field $W$ as the difference of two connections, the one being the background connection, then  $W^{k} :=g^{mn}(\Gamma^{k}_{mn}-\tilde{\Gamma}^{k}_{mn})$ is a tensor.
 Let us express (\ref{Veq})  in a different form as 
 \begin{equation}
	V_{j}=\frac{1}{2}\nabla^{k}\left [D(Ric_{g})(h)\right]_{jk}-\Delta W_{j}.
	\end{equation}
By the same reasoning, we can write
\begin{equation}
\left[D(Ric_{g})(h)\right]_{jk}=-\frac{1}{2} \Delta h_{jk} +\frac{1}{2}(\nabla_{j}\tilde V_{k}+\nabla_{k}\tilde V_{j}),
\end{equation}
where 
\begin{equation}
\tilde V_{j}=g_{jk}g^{mn}[D(\Gamma_{g})(h)]^{k}_{mn}.
\end{equation}
Again, we may define $W^{k}$ as the difference between two connections to be a vector field. What we finally have is
\begin{eqnarray}
V_{j}&=&\frac{1}{4}\nabla^{k} (- \Delta h_{jk} +\nabla_{j}W_{k}+\nabla_{k}W_{j})-\Delta W_{j} \nonumber \\
&=& -\frac{1}{4}\nabla^{k}\Delta h_{jk}+\frac{1}{4}\nabla^{k}\nabla_{j}W_{k}-\frac{3}{4}\Delta W_{j}.
\end{eqnarray}
We once more write the linearized $K$-tensor but call back +$W$ terms to the equation
 \begin{eqnarray}
	\left[D(K_{g}(h))\right]_{ij}
	&=&\Delta^{2}h_{ij} + \frac{1}{2} g_{ij}(\Delta\nabla_{k}\nabla_{l}h^{kl}-\Delta^{2}h)-\frac{1}{4}\left(\nabla_{i}\nabla^{k}\Delta h_{jk}+\nabla_{j}\nabla^{k}\Delta h_{ik}\right) \nonumber \\
	&+&\frac{1}{4} \left (-\nabla_{i}\nabla_{j}	\Delta h+2\nabla_{i}\nabla_{j}\nabla_{k}\nabla_{l}h^{kl}\right)-\frac{3}{4}\Delta (\nabla_{i} W_{j}+\nabla_{j} W_{i}) 	.  
	\end{eqnarray}
 We expect to cancel the contributions from $-\Delta W$ terms which are of fourth-order derivatives of the metric, so that linearized $K$ tensor is an elliptic operator, more properly a bi-laplacian. As deTurck did, we can define a one-parameter family of diffeomorphisms $\phi_{t}$  induced by the vector field $W(t)$ which is the generator of the integral curves induced by $\phi_{t}$ since we consider compact manifolds; and for compact manifolds vector fields and integral curves are complete. Then for the contributions from $W$ terms, we can write $\nabla_{i}W_{j}+\nabla_{j}W_{i}=(\mathcal{L}_{W}g)_{ij}$ that is the Lie derivative of the metric $g(t)$ with respect to the vector field $W(t)$. By this construction, the linearized tensor reads
 \begin{eqnarray}
	\left[D(K_{g}(h))\right]_{ij}&=&\Delta^{2}h_{ij} + \frac{1}{2} g_{ij}(\Delta\nabla_{k}\nabla_{l}h^{kl}-\Delta^{2}h)-\frac{1}{4}(\nabla_{i}\nabla^{k}\Delta h_{jk}+\nabla_{j}\nabla^{k}\Delta h_{ik}) \nonumber \\
	&+&\frac{1}{4} \left(-\nabla_{i}\nabla_{j}	\Delta h+2\nabla_{i}\nabla_{j}\nabla_{k}\nabla_{l}h^{kl}\right)-\frac{3}{4}\Delta(\mathcal{L}_{W}g)_{ij}. 
	\end{eqnarray}
 To break the diffeomorphisms we may use two constructions. The first construction is to use the Lie derivative of the connections, we define the connection as being the difference of two connections again so that it is a tensor also. Now we have
\begin{equation}
\mathcal{L}_{W}\Gamma^{k}_{ij}= \nabla_{i}\nabla_{j}W^{k}+R^{\phantom{ilj}k}_{ilj}W^{l}.
\end{equation}
On the other hand, since the Lie derivative of a vector field with respect to itself vanishes we can write
\begin{eqnarray}
\mathcal{L}_{W} \left(g^{km}\Gamma^{l}_{km} \right)&=&\mathcal{L}_{W}(g^{km})\Gamma^{l}_{km}+g^{km}\mathcal{L}_{W}\Gamma^{l}_{km}=0 \nonumber \\
&=&-2\Gamma^{l}_{km}\nabla^{k}W^{m}+\Delta W^{l}+R^{\phantom{n}l}_{n}W^{n}, \nonumber \\
\Delta W^{l}&=&2\Gamma^{l}_{km}\nabla^{k}W^{m}-R^{\phantom{n}l}_{n}W^{n}, \nonumber \\
\Delta \nabla_{i} W_{j}&=&\nabla_{i}g_{jl} \Delta W^{l}=\nabla_{i}(2g_{jl}\Gamma^{l}_{km}\nabla^{k}W^{m})- \nabla_{i}(R_{nj}W^{n}).		
\end{eqnarray}
For the principal symbol, we need only the highest-order derivatives, so that no contribution comes from the lower-order terms. Then the principal symbol becomes
\begin{eqnarray}
	\hat{\sigma}\left[D(K)\right](\zeta)(h)_{ij}&=& |\zeta|^{2} |\zeta|^{2} h_{ij} +\frac{|\zeta|^{2}}{2}g_{ij} \left(\zeta_{k}\zeta_{l}h^{kl}-  |\zeta|^{2}h \right) \nonumber \\
						&-&\frac{|\zeta|^{2}}{4} \left(\zeta_{i}\zeta^{k}h_{jk}+\zeta_{j}\zeta^{k} h_{ik}\right)+\frac{ \zeta_{i}\zeta_{j}}{4}\left (-|\zeta|^{2}h+2\zeta_{k}\zeta_{l}h^{kl}\right).
\end{eqnarray}
Choosing the problematic gauge mode $h_{ij}=\zeta_{i}\zeta_{j}$, we have 
	\begin{equation}
	\langle \hat{\sigma}[D(K)](\zeta)(h),h\rangle=\frac{3}{4} |\zeta |^{2} |\zeta|^{2}  |\zeta |^{2} |\zeta |^{2} \neq 0.
	\end{equation}
Therefore $D(K)$ is an elliptic operator, and we conclude that the unspecified constant $\alpha$ in the flow equation must be chosen to be negative so that the flow is diffusive as expected, and by rescaling $t$ we may set  $\alpha= -1$. [Note that if we had taken the form of the $K$-tensor as given in (\ref{K-tensor}), then we should have set $\alpha =+1$.]

Another construction of the short-time existence and uniqueness was given in \cite{Jeff};\footnote{We were not aware of this work when we submitted the first version of our paper to the arXiv; we duly thank J.D. Streets for letting us know about his work.} and we simply follow the results therein for the benefit of the reader. Consider two Riemannian metrics $g(t)$ and $\bar{g}(t)$ with $\phi_{t}^{*}g(t)=\bar{g}(t)$, and $\phi_{t}$ a one-parameter family of diffeomorphisms induced by the vector field $V(t)=-\frac{3}{4}\Delta W(t)$ so that
\begin{eqnarray}
\frac{\partial \phi_{t}(x)}{\partial t}&=&-V\Big(\phi_{t}(x),g(t),t\Big), \nonumber \\
\phi_{t}(0)&=& id_{M}  \label{diffV} \quad .
\end{eqnarray}
The vector field $W^{k}$ can be written as 
\begin{equation}
-W^{k}=\Big( (\phi_{t}^{-1}) ^{*}\bar{g}    \Big)^{mn}(-\Gamma^{k}_{mn}+\tilde{\Gamma}^{k}_{mn})=(\Delta_{\bar{g}(t),\tilde{g}}\phi_{t})^{k} ,
\end{equation}
where $\tilde{g}$ is a fixed background metric on the manifold, and $\Delta_{\bar{g}(t),\tilde{g}}\phi_{t}$ is the harmonic map Laplacian  with respect to the domain metric $\bar{g}$ and the codomain metric $\tilde{g}$ which can be computed by using the equation for the Christoffel symbol of the pullback metric  
\begin{equation}
-g^{\alpha \beta}\Gamma(g)^{\gamma}_{\alpha \beta}=(\phi^{*}g)^{ij}\Big( \frac{\partial^{2}\phi^{\gamma}}{\partial x^{i}\partial x^{j}}    -\Gamma^{k}_{ij} (\phi^{*}g) \frac{\partial \phi^{\gamma}}{\partial x^{k}} \Big) , 
\end{equation}
as given in  \cite{CK}. Since the vector field $V$ is a Laplacian of the vector field $W$, we need to relate the  Laplacians with respect to the two metrics $g$ and $\bar{g}$ acting on the one-form $W$, and one can show that the two Laplacians differ  by third-order derivatives of the map $\phi$ at the highest order
\begin{equation}
\Delta_{g} W = \Delta_{\bar{g}} W + {\mathcal{O}}(\partial^{3} \phi) . \label{lapl}
\end{equation}
By using equation (\ref{lapl}), we may write an ODE for the diffeomorphisms in terms of $\bar{g}(t)$ as
\begin{eqnarray}
\frac{\partial \phi_{t}(x)}{\partial t}&=&L(\phi_{t}(x),\bar{g}(t),t), \nonumber \\
\phi_{t}(0)&=& id_{M}  \label{diffop}, \quad 
\end{eqnarray}
where $L$ is a fourth-order parabolic operator, a bi-laplacian for $\phi_{t}$ with solutions existing at least for a short time if $\bar{g}(t)$ exists which is the solution of the initial value problem 
\begin{equation}
\partial_{t}\bar{g}(t)=-K(\bar{g}(t)) \qquad \bar{g}(0)=g_{0}. \label{IVP} 
\end{equation}
If $\bar{g}(t)$ exists, then  $g(t)=(\phi_{t})_{*}\bar{g}(t)$ is a solution of the gauge-fixed flow 
\begin{equation}
\partial_{t}g(t)=-K(g(t))+ \mathcal{L}_{V}g(t) \qquad g(0)=g_{0} , \label{gauged} 
\end{equation}
which can be computed as follows
\begin{eqnarray}
\frac{\partial}{\partial_{t}}\bar{g}(t)&=&\frac{\partial}{\partial_{t}}(\phi_{t}^{*}g(t)) 
=\frac{\partial}{\partial_{s}} \bigg|_{s=0}    \phi_{t+s}^{*}g(t+s) =     \phi_{t}^{*}     \frac{\partial}{\partial_{s}} \bigg|_{s=0} g(t+s) +  \frac{\partial}{\partial_{s}} \bigg|_{s=0} \phi_{t+s}^{*}g(t) \nonumber \\
 &=&  \phi_{t}^{*} (  -K[g(t)])+ \phi_{t}^{*}  \mathcal{L}_{V(t)}g(t)+\frac{\partial}{\partial_{s}}\bigg|_{s=0}\biggl( \Bigl( \phi_{t}^{-1}\circ \phi_{t+s} \Bigr)^{*}\phi_{t}^{*}g(t) \biggr) \nonumber \\
 &=&-K[\phi_{t}^{*}g(t)]+\phi_{t}^{*}\Bigl( \mathcal{L}_{V(t)}g(t)\Bigr) -      \mathcal{L}_{(\phi_{t}^{-1})_{*}V(t)} \phi_{t}^{*}g(t) \nonumber \\
&=&-K[\bar{g}(t)] .
\label{Diffeo}
\end{eqnarray}
The gauge-fixed flow (\ref{gauged}) is elliptic, it has a bi-laplacian symbol, so will have solutions at least for a short time. Then $\bar{g}(t)$ exits and solves the $K$-flow (\ref{IVP}). 

Let us assume that  $\bar{g}_{1}(t)$ and  $\bar{g}_{2}(t)$ are two solutions of the $K$-flow (\ref{IVP}) on a common time interval with $\phi_{1}(t)$ and $\phi_{2}(t)$ being the solutions of the flow (\ref{diffop}) with respect to the domain metrics $\bar{g}_{1}(t)$ and  $\bar{g}_{2}(t)$ respectively, and the codomain metric $\tilde{g}$. Then
\begin{equation}
g_{i}(t)=((\phi_{i})_{t})_{*}\bar{g}_{i}(t), \qquad i=1,2, \quad \nonumber
\end{equation}
are solutions of the gauge-fixed flow (\ref{gauged}).  Since $g_{1}(0)=g_{2}(0)$ and (\ref{gauged}) have unique solutions, then $g_{1}(t)=g_{2}(t)$ as long as both exist. The flow of $\phi (t)$ can also be written in terms of $g(t)$ as in (\ref{diffV}). Since $g_{1}(t)=g_{2}(t)$ on a common time interval, $\phi_{1}(t)=\phi_{2}(t)$ on a common time interval as well which then implies $\bar{g}_{1}(t)=\bar{g}_{2}(t)$  and the uniqueness of  $\bar{g}_{i}(t)$ is inherited by the equivalence of (\ref{diffV}) and (\ref{diffop}).

\section{Gradient flow and Entropy Formulation of the flow }
As alluded to in the Introduction section, the flow is related to the quadratic part of the new massive gravity theory. Let us now show that the action of the theory is non-decreasing along the flow, hence it can be considered as an entropy. 
\begin{equation}
S=-\int_{\mathcal{M}} d^{3} x \sqrt{g} \left ( \dfrac{3}{8} R^{2}- R_{i j}R^{i j}\right )=-\int_{\mathcal{M}} d^{3} x \sqrt{g} K,
\label{action}
\end{equation}
where ${\mathcal{M}}$ is a closed 3-manifold, so we can drop the boundary terms after integration by parts.
By taking its variation with respect to the metric we obtain
\begin{equation}
\delta_g S=-\int_{\mathcal{M}} d^{3} x \sqrt{g}  ( K^{ij})\delta g_{ij}  
=\int_{\mathcal{M}} d^{3} x \sqrt{g} K^{ij} K_{ij}, 
\end{equation}
where we used $\delta g_{ij}=-K_{ij}$. Thus the action is a non-decreasing functional of time, as we would require from an entropy functional. If we were to consider infinitesimal diffeomorphisms and scalings $\delta g_{ij}=\nabla_{i}X_{j}+\nabla_{j}X_{i}+\lambda(x) g_{ij}$, then  we would find 
\begin{equation}
\delta S=-\int_{\mathcal{M}}  d^{3} x \sqrt{g} \lambda K, 
\end{equation}
which shows that diffeomorphisms are symmetries but the scalings are not. See Appendix D where the transformation properties of the $K$ tensor under conformal scalings were worked out.  Note also that we have used the Bianchi Identity $\nabla_i K^{i j}=0$ which is valid for all smooth metrics. 

The functional (\ref{action}) has been thoroughly studied by  Gursky and Viaclovsky  \cite{Viac}, and here we follow their main results for definitions and theorems. 

Let $\mathcal{H}$ denote the set of smooth Riemannian metrics on a compact manifold  $\mathcal{M}$ and let $\mathcal{H}_{1}$ be the set of smooth Riemannian metrics that keep the volume constant, set to 1.  A real-valued functional on $\mathcal{H}_{1}$, is called a Riemannian functional $F$, if it is invariant under the action of the diffeomorphism group, that is $F(g)=F(\phi^{*}g)$ for $ g\in C^{\infty}(S_{2}T^{*}M)$ and $\phi$ is an element of the diffeomorphism group. Every Riemannian functional (constructed from a polynomial of the Riemann curvature tensor $\text{Riem}$) is differentiable. Then the differential of $F$ is defined as
\begin{equation}
\left. \frac{d}{dt}\right \vert_{t=0}F(g+th)=F^{'}_{g}\cdot h=\int_\mathcal{M} g(h,\nabla F)d\mu,
  \end{equation}
where $\nabla F$ is called the gradient of $F$ and it is a symmetric two tensor and a polynomial of Riemann curvature tensor $\text{Riem}$ and its first two covariant derivatives $\nabla \nabla \text{Riem}$. The gradient of a Riemannian functional is covariantly conserved; it defines a vector field on $\mathcal{H}_{1}$ and this vector field satisfies an integrability condition. For the action (\ref{action}), the gradient is the tensor $K_{ij}$.

\subsection{Results of  Gursky and Viaclovsky  \cite{Viac}}

In \cite{Viac}, the authors consider the negative of the action (\ref{action}), and call it $F_{2}$ so $F_{2}=-S=\int d\mu K$. We state here their results in the context of compact three-dimensional manifolds, with the restriction that the volume is kept fixed, that is we consider $S \mid_{\mathcal{H}_{1}}$.

If $ K>0$, then the sectional curvatures of $g$ at a point are all positive or all negative. Now consider two cases.
\begin{enumerate}
\item  A metric $g$ with $F_{2}>0 \quad (S<0)$ is critical for for $F_{2} \mid_{\mathcal{H}_{1}} (S \mid_{\mathcal{H}_{1}})$ if and only if $g$ has constant sectional curvatures, strictly positive or negative. This condition is called {\it ellipticity}, and respectively critical points are called elliptic critical points. Then elliptic critical points of $S$ are of constant curvature. 

\item  A metric $g$ with $F_{2}=0 \quad (S=0)$ is critical for $F_{2} \mid_{\mathcal{H}_{1}} (S \mid_{\mathcal{H}_{1}})$ if and only if $g$ has constant sectional curvature. This condition is called {\it degenerate ellipticity}. In this case, it is possible that the curvature changes sign. 

\end{enumerate}

For both cases, a metric is critical if and only if it has constant sectional curvatures, so constant curvature. A Riemannian metric on a 3-dimensional manifold is Einstein if and only if 
it has constant sectional, hence critical metrics are Einstein metrics.\footnote{ In 3 dimensions, the relation between the Riemann curvature tensor and the Einstein tensor is as $R_{\mu \alpha\nu\beta} =\frac{1}{4}\epsilon_{\mu\alpha\sigma}\epsilon_{\nu \beta \rho}G_{\sigma\rho}$.}; Moreover, if a 3-dimensional manifold admits a constant sectional curvature metric, then its universal cover is diffeomorphic to $\mathbb{R}^{3}$ or $S^{3}$.
 Constant curvature metrics put a condition on the trace $K$. Since, for a constant curvature metric $\text{Ric}=\frac{1}{3}R g$, then $K=\frac{1}{24}R^{2}\geq 0 \quad (-K\leq 0)$. If $g$ is critical for $S \mid_{\mathcal{H}_{1}}$, then $K$ is a constant.
\begin{itemize}

\item If  $S<0$, then  $K>0$, or $(-K<0)$. The scalar curvature $R$ can be positive or negative, but constant.  

\item  If $S=0$, then $K=0$ and $R\leq 0$ and $g$ is conformally flat. 
\end{itemize}

Finally, for $F_{2}<0 \quad (S>0)$, there are critical metrics that are called non-elliptic critical points. In \cite{Viac}, one such example on $S^{3}$ is given, which is the orthonormal frame metric.  

In conclusion for $S\leq 0$, such that  $K\geq 0$, the elliptic critical points of $S$ are either constant curvature metrics, positive or negative, or just conformally flat metrics. Thus under the $K$-Flow, any topological constant curvature manifold will approach the geometric manifold with the standard metric. For $S\geq 0$, it is possible to find critical metrics but a complete classification is not possible. For example for $\mathbb{R}^{3}$ we know that $K=0$, and under the flow the geometry does not change, it is still flat, but it is trivial.

\section{Some solutions of the $K$-flow}

\subsection{Einstein Solutions}

In this part of our work, we discuss some specific solutions to the flow. We start with an Einstein manifold, a smooth manifold that admits an Einstein metric $g_{ij}$ such that the Ricci tensor is simply $R_{ij}=2 \Lambda g_{ij}$ where $\Lambda$ is just a constant. If initially $R_{ij}(x,0)= 2\Lambda  g_{ij}(x,0)$ for all $x
\in M$, then
\begin{equation}
R(x,0)=6\Lambda, \quad S_{ij}(x,0)=\frac{1}{2}\Lambda g_{ij}(x,0), \quad K_{ij}(x,0)=\frac{1}{2} \Lambda^{2}  g_{ij}(x,0).
\end{equation}
Let us consider the 3-cases separately.
\begin{enumerate}
 
\item  Consider the $\Lambda=0$ case, that is initially the metric is Ricci-flat (or a Milnor frame metric). Then $K_{ij}(x,0)=0$ and so we have a stationary solution of the flow. 

\item  For $\Lambda >0$, initially the scalar curvature is positive, and the metric will shrink under the flow. To see this, we may set 
\begin{equation}
g_{ij}(x,t)=a^{2}(t)g_{ij}(x,0)
\end{equation}
for some smooth function $a(t)$. Then the relevant tensors  and geometric objects read
\begin{eqnarray}
&g^{ij}(x,t)&=a^{-2}(t)g^{ij}(x,0), \qquad \Gamma^{k}_{ij}(x,t)=\Gamma^{k}_{ij}(x,0), \qquad R_{\phantom{k}ilj}^{k}(x,t)=R_{\phantom{k}ilj}^{k}(x,0), \nonumber \\
 &R_{ij}(x,t)&=R_{ij}(x,0), \qquad R(x,t)=a^{-2}(t)R(x,0), \qquad S_{ij}(x,t)=S_{ij}(x,0), \nonumber \\
 &K_{ij}(x,t)&=a^{-2}(t)K_{ij}(x,0)=\frac{1}{2} \Lambda^{2} a^{-2}  g_{ij}(x,0).
\end{eqnarray}
Therefore, the flow equation becomes
\begin{equation}
\partial_{t}(a^{2}(t)g_{ij}(x,0))=-\frac{1}{2} a^{-2}\Lambda^{2}  g_{ij}(x,0),
\end{equation}
with the solution 
\begin{equation}
a^{4}(t)=a^{4}_{0}-\Lambda^{2} t.
\end{equation}
The metric $g_{ij}(x,t)$ will shrink to a point in a finite time $t=a_0^{4} / \Lambda^{2}$, the rate of decay depends on the parameters $a_{0}$ and $\Lambda$. But eventually, the scalar curvature becomes infinite. 
\item  For $\Lambda <0$, that is for an initially negative scalar curvature Einstein metric, the flow equation would be exactly as in the second case above, since $K$ is proportional to $\Lambda^{2}$ term.  Thus, the metric $g_{ij}(x,t)$ shrinks and the scalar curvature becomes negatively infinite in a finite amount of time.

\end{enumerate}

\subsection{Ancient Solutions}
Consider the three dimensional round sphere with the standard metric  $g_{\text{can}}$, and one parameter family of conformally equivalent metrics $g(t)=r^{2}(t)g_{can}$, then $g(t)$ is a solution of the flow if 
\begin{equation}
\partial_{t} g(t)=2r(t) \dfrac{dr(t)}{dt} g_{can}=- K[g(t)]=- \dfrac{r^{-2}}{2}(t) g_{can},
\end{equation}
with the solution
\begin{equation}
r(t)=(r_{0}^{4}-t)^{1/4}.
\end{equation}
Hence this solution exists for $(S^{3},g(t))$ in the time interval $(-\infty,T)$ where $T$ is the singularity time $T=r_{0}^{4}$, hence an ancient solution. This is an example of the second case above.

\subsection{The Model Geometries Under $K$-Flow}
In this part, we will study the properties of the geometries under the flow. Since the $K$-tensor is not traceless we can not set the volume to a constant value say 1, in other words, the volume is not preserved under the flow. We could of course write a traceless form of the tensor and work with this traceless tensor in which case we could set the volume to 1. However in this case the action/entropy would vanish. Thus we do not consider a normalized flow. We start with the {\it Non-Bianchi classes} for which we use the coordinate frame and then move to  {\it Bianchi classes} and study them in the orthonormal frame. 
\vspace{0.5cm}
\subsubsection{Non-Bianchi Classes}

Now we work on the so-called Non-Bianchi Classes whose metrics in the coordinate basis are simple enough to yield exact integration of the flow equations. Thus we study their flows on a coordinate basis. 

\vspace{0.5cm}
{\bf 1: The geometry of $\mathbb{H}^{3}$ with the metric $g=a(t) g_{\mathbb{H}^{3}}$}	
\vspace{0.5cm}

The metrics of the geometries in this class are of the above form where $g_{\mathbb{H}^{3}}$ is the metric of $\mathbb{H}^{3}$ with a smooth function $a(t)$. In the basis $\left\lbrace dx, dy, dz \right\rbrace $  the components of the Ricci tensor and the curvature scalar can be computed directly for the metric
\begin{equation}
{\bf g}=\frac{a(t)}{x^{2}}(dx^{2}+   dy^{2}  +dz^{2}).
\end{equation}
One has $\partial_x g_{ij}=-2a\delta_{ij}/x^{3}$ as the only non-zero derivatives of the components of the metric. Therefore the non-zero Christoffel symbols are
\begin{equation}
\Gamma^{x}_{xx}=\Gamma^{y}_{xy}=\Gamma^{z}_{xz}=-\frac{1}{x}, \qquad \Gamma^{x}_{yy}=\Gamma^{x}_{zz}=\frac{1}{x},
\end{equation}which lead to 
\begin{equation}
R_{ij}= -\frac{2}{x^{2}}\delta_{i j}, \qquad R= -\frac{6}{a}, \qquad K_{ij}=\frac{1}{2ax^{2}}\delta_{i j},\qquad K=\frac{3}{2a^{2}},
\end{equation}
the flow yields
\begin{eqnarray}	
		a^{2}(t)&=& a^{2}_{0} -t.
	\end{eqnarray}
	In this case, all three dimensions shrink with $t^{1/2}$ and the scalar curvature becomes negative infinite while the trace $K$ also becomes positive infinite in a finite amount of time.
	
\vspace{0.5cm}
{\bf 2: The geometry of $\mathbb{H}^{2}\times \mathbb{R}$ with the metric $g=a(t)g_{\mathbb{R}} + b(t) g_{\mathbb{H}^{2}}$}
\vspace{0.5cm}

The metrics of the geometries in this class are of the above form where $g_{\mathbb{R}}$ is the metric of $\mathbb{R}$, and $g_{\mathbb{H}^{2}}$ is the metric of the hyperbolic plane with smooth functions $ a(t)$  and $ b(t)$ respectively. 

In the  basis $\left\lbrace dx,d\theta,\sinh\theta d\phi\right\rbrace, $ the metric reads
\begin{equation}
{\bf g}=a(t)dx^{2}+b(t)(d\theta^{2}+\sinh^{2}\theta d^{2}\phi),
\end{equation}
from which one computes $\partial_\theta g_{\phi\phi}=2a \sinh\theta \cosh\theta $ which is the only non-zero derivative of the components of the metric. 
The non-zero Christoffel symbols are
\begin{equation}
\Gamma^{\theta}_{\phi \phi}=-\sinh\theta \cosh\theta, \qquad \Gamma^{\phi}_{\phi \theta}=\Gamma^{\phi}_{\theta \phi }=\dfrac{\cosh\theta}{\sinh\theta},
\end{equation}
yielding
\begin{equation}
R_{xx}=0, \qquad R_{\theta \theta}=-1, \qquad R_{\phi \phi}=-\sinh^{2}\theta, \qquad R=-\dfrac{2}{b},
\end{equation}
and
\begin{equation}
K_{11}= \frac{a}{2b^{2}}, \qquad
K_{22}= -\frac{1}{2b}, \qquad
K_{33}= -\frac{\sinh^{2}\theta}{2b} \qquad K= -\frac{1}{2b^{2}}.
\end{equation}
Finally, the flow equations read :
\begin{equation}
\partial_{t}a=-\frac{a}{2b^{2}},  \qquad \partial_{t}b=\frac{1}{2b},
\end{equation}
which can be directly integrated
\begin{equation}
b^{2}(t)=b_{0}^{2}+t,  \qquad a(t)=\frac{a_{0}b_{0}}{\sqrt{b_{0}^{2}+t}}.
\end{equation}
From these solutions, one can see that, for large times, while the hyperbolic plane expands linearly with  $t^{1/2}$,  $\mathbb{R}$ shrinks with $t^{-1/2}$ giving a {\it pancake degeneracy} to the geometry. The scalar curvature vanishes with $t^{-1/2}$. The trace $K$ also vanishes. 

\vspace{0.5cm}	
{\bf 3: The geometry of $S^{2}\times \mathbb{R}$ with the metric $g=a(t)g_{\mathbb{R}} + b(t) g_{S^{2}}$}
\vspace{0.5cm}

The metrics of the geometries in this class are of the above form where $g_{\mathbb{R}}$ is the metric of $\mathbb{R}$  and $g_{S^{2}}$ is the metric of the two-sphere with smooth functions $ a(t)$  and $ b(t)$ respectively. In the basis, $\left\lbrace dx,d\theta,\sin\theta d\phi\right\rbrace $, the metric is 
\begin{equation}
{\bf g}=a(t)dx^{2}+b(t) (d\theta^{2}+\sin^{2}\theta d^{2}\phi),
\end{equation} 
and $\partial_\theta g_{\phi\phi}=2b \sin\theta \cos\theta $ is the only non-zero derivative of the components of the metric. One has
\begin{equation}
\Gamma^{\theta}_{\phi \phi}=-\sin\theta \cos\theta, \qquad \Gamma^{\phi}_{\phi \theta}=\Gamma^{\phi}_{\theta \phi }=\frac{\cos\theta}{\sin\theta},
\end{equation}
yielding
\begin{equation}
R_{xx}=0, \qquad R_{\theta \theta}=1, \qquad R_{\phi \phi}=\sin^{2}\theta, \qquad 
R=\dfrac{2}{b}. 
\end{equation}
Then we find
\begin{equation}
K_{11}= \frac{a}{2b^{2}}, \qquad K_{22}= -\frac{1}{2b}, \qquad K_{33}= -\frac{\sin^{2}\theta}{2b}, \qquad K= \frac{-1}{2b^{2}}.
\end{equation}
The flow equations yield the same equations we obtained for the previous $\mathbb{H}^{2}\times \mathbb{R}$ case. Thus for large times, the two-sphere expands linearly with  $t^{1/2}$, and  $\mathbb{R}$ shrinks with $t^{-1/2}$. The scalar curvature vanishes as $t^{-1/2}$. It is interesting that  $S^{3}$ and $\mathbb{H}^{3}$ shrink to a point, but  $S^{2}$ and $\mathbb{H}^{2}$, on the other hand, expand, in the $K$-flow in contrast to the case in Ricci flow. Ricci flow gives shrinking solution for all $S^{n}$ and expanding solution for all $\mathbb{H}^{n}$, independent of the number of dimensions. 
\begin{figure}[h]
    \centering
    \includegraphics[width=110mm,scale=0.5]{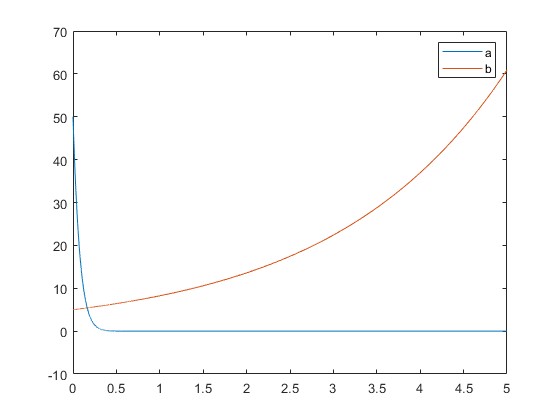}
    \caption{The behavior of both $\mathbb{H}^{2}\times \mathbb{R}$ and $S^{2}\times \mathbb{R}$ geometries under the flow is plotted. The vertical axis represents the dimensions $a,b$. The horizontal axis represents the parameter of the flow $t$.}
    \label{Fig.H(2)R}
\end{figure}

\subsubsection{The Bianchi class geometries: Orthonormal Frame Calculations}	
We now consider the geometries known as the Bianchi classes and describe their behavior under the flow. We have to find the $K$-tensor (or more properly the $K$ one-form) in the orthonormal frame. This is a straightforward but rather tedious computation. Therefore we delegate the bulk of this computation to Appendix A and introduce here the essentials needed to write down the coupled non-linear ODEs coming from the flow.

We now introduce the Milnor frame and compute the relevant geometrical quantities there. 
Let ${\bf g}$ be a left-invariant metric on the unimodular, simply connected Lie group \emph{G}, then there exists a left-invariant orthogonal frame   $F=\left\lbrace F_{i}\right\rbrace $, the Milnor frame, such that :
\begin{align}
[F_{i},F_{j}]=C_{ij}~^{k}F_{k}, 
\end{align}
where $C_{ij}^{\phantom{ij}k}$ are the structure constants uniquely defined for each Lie algebra of the Lie groups that we shall consider. Let $ \nabla$ be a uniquely defined connection satisfying the identities: 
\begin{align}
&\nabla_{X}Y-\nabla_{Y}X=[X,Y], \nonumber \\
&\left\langle \nabla_{X}Y,Z\right\rangle +\left\langle Y,\nabla_{X}Z\right\rangle =0, \nonumber \\
&\left\langle \nabla_{X}Y,Z\right\rangle=\dfrac{1}{2} \Big ( \left\langle [X,Y],Z\right\rangle  - \left\langle [Y,Z],X\right\rangle +\left\langle [Z,X],Y\right\rangle \Big ), 
\end{align}
where $\langle Y,Z\rangle$ denotes the inner product and $X,Y,Z$ are vector fields. In an orthonormal coframe $\left\lbrace E_{i}\right\rbrace $ , where $[E_{i},E_{j}]=D_{ij}^{\phantom{ij}k}E_{k}$, one has  
\begin{align}
\left\langle \nabla_{E_{i}}E_{j},E_{k}\right\rangle=\left\langle [E_{i},E_{j}],E_{k}\right\rangle =\dfrac{1}{2}\left ( D_{ij}~^{k}-D_{jk}~^{i}+D_{ki}~^{j}\right ).
\end{align}
Hence $\nabla_{E_{i}}E_{j}=\sum_{k} \dfrac{1}{2}\left ( D_{ij}^{\phantom{ij}k}-D_{jk}^{\phantom{jk}i}+D_{ki}^{\phantom{ki}j}\right)E_{k}$. $\nabla_{E_{i}}E_{j}= -\nabla_{E_{j}}E_{i}$ and $\nabla_{E_{i}}E_{i}=0$ which follows also from the properties of the structure constants. Conventionally one writes the algebra explicitly as
\begin{align}
[F_{1},F_{2}]=\nu F_{3}, \quad [F_{2},F_{3}]=\lambda F_{1} , \quad [F_{3},F_{1}]=\mu F_{2},
\end{align}
where $\mu ,\nu, \lambda \in \left\lbrace -1,0,1\right\rbrace $ are the structure constants. In this frame, then the metric ${\bf g}$ can be written as :
\begin{align}
{\bf g}=a(t)\omega^{1}\otimes\omega^{1}+b(t)\omega^{2}\otimes\omega^{2}+c(t)\omega^{3}\otimes\omega^{3},
\end{align}
where $\left\lbrace \omega^{1},\omega^{2},\omega^{3}\right\rbrace $ are the dual basis one-forms to the orthogonal basis $\left\lbrace F_{1},F_{2},F_{3}\right\rbrace $ , that is $\omega^{i}(F_{j})=\delta^{i}~_{j}$ . Let then the  dual basis one-forms to the orthonormal basis  $\left\lbrace E_{1},E_{2},E_{3}\right\rbrace $ be $\left\lbrace e^{1},e^{2},e^{3}\right\rbrace $. In this orthonormal coframe, the metric takes the form : 
\begin{align}
{\bf g}=e^{1}\otimes e^{1}+e^{2}\otimes e^{2}+e^{3}\otimes e^{3},
\end{align}
where the relations between these bases are :
\begin{align}
e^{1}=\sqrt{a}\omega^{1}, \qquad e^{2}=\sqrt{b}\omega^{2}, \qquad e^{3}=\sqrt{c}\omega^{3}, 
\end{align} 
and since $e^{i}(E_{j})=\delta^{i}~_{j}$ we have,
\begin{align}
E_{1}=\dfrac{F_{1}}{\sqrt{a}}, \qquad E_{2}=\dfrac{F_{2}}{\sqrt{b}}, \qquad E_{3}=\dfrac{F_{3}}{\sqrt{c}} .
\label{EFmilnor}
\end{align}

As it is given in  Appendix A, one finds the $K$ one-forms in the orthonormal frame as 
	\begin{eqnarray}
		K^{1}&=& \frac{-1}{32 a^2 b^2 c^2}\Bigg (-105a^4\lambda^4+60a^3b\lambda^3\mu+60a^3c\lambda^3\nu+2a^2b^2\lambda^2\mu^2-20a^2bc\lambda^2\mu\nu \nonumber \\ 
 &+&2a^2c^2\lambda^2\nu^2 -20ab^3\lambda\mu^3+20ab^2c\lambda\mu^2\nu+20abc^2\lambda\mu\nu^2-20ac^3\lambda\nu^3+63b^4\mu^4\nonumber  \\ 
&-&60b^3c\mu^3\nu-6b^2c^2\mu^2\nu^2-60bc^3\mu\nu^3+63c^4\nu^4 \Bigg ) e^{1}, \\
		K^{2}&=& \frac{-1}{32 a^2 b^2 c^2} \Bigg (63a^4 \lambda^4-20a^3 b \lambda^3 \mu -60a^3 c \lambda^3 \nu + 2a^2 b^2 \lambda^2 \mu^2 +20a^2 b c \lambda^2 \mu \nu  \nonumber \\ 
 &-&6a^2 c^2 \lambda^2 \nu^2+60ab^3 \lambda \mu^3 -20ab^2 c \lambda \mu^2 \nu +20abc^2 \lambda \mu \nu^2 -60ac^3 \lambda \nu^3-105b^4 \mu^4  \nonumber \\
&+&60b^3 c \mu^3 \nu +2b^2 c^2 \mu^2 \nu^2-20bc^3\mu\nu^3+63c^4\nu^4 \Bigg )e^{2},\\
		K^{3}&=& \frac{-1}{32 a^2 b^2 c^2}\Bigg(63a^4\lambda^4-60a^3b\lambda^3\mu-20a^3c\lambda^3\nu-6a^2b^2\lambda^2\mu^2+20a^2bc\lambda^2\mu\nu \nonumber \\  
&+&2a^2c^2\lambda^2\nu^2-60ab^3\lambda\mu^3+20ab^2c\lambda\mu^2\nu-20abc^2\lambda\mu\nu^2+60ac^3\lambda\nu^3 \nonumber \\  
&+&63b^4\mu^4-20b^3c\mu^3\nu+2b^2c^2\mu^2\nu^2+60bc^3\mu\nu^3-105c^4\nu^4\Bigg) e^{3}.
	\end{eqnarray}
Then using the flow equation $\partial_t e^a= - K^a$, we have the following coupled nonlinear ODEs:
	\begin{eqnarray}
		\frac{da}{dt}&=& \frac{a}{16 a^2 b^2 c^2} \Bigg(-105a^4\lambda^4+60a^3b\lambda^3\mu+60a^3c\lambda^3\nu+2a^2b^2\lambda^2\mu^2-20a^2bc\lambda^2\mu\nu \nonumber \\  
&+&2a^2c^2\lambda^2\nu^2-20ab^3\lambda\mu^3+20ab^2c\lambda\mu^2\nu+20abc^2\lambda\mu\nu^2-20ac^3\lambda\nu^3+63b^4\mu^4 \nonumber \\
&-&60b^3c\mu^3\nu-6b^2c^2\mu^2\nu^2-60bc^3\mu\nu^3+63c^4\nu^4 \Bigg),  \\
		\frac{db}{dt}&=& \frac{b}{16 a^2 b^2 c^2}\Bigg(63a^4\lambda^4-20a^3b\lambda^3\mu-60a^3c\lambda^3\nu+2a^2b^2\lambda^2\mu^2+20a^2bc\lambda^2\mu\nu \nonumber  \\ 
&-&6a^2c^2\lambda^2\nu^2+60ab^3\lambda\mu^3-20ab^2c\lambda\mu^2\nu +20abc^2\lambda\mu\nu^2-60ac^3\lambda\nu^3-105b^4\mu^4 \nonumber \\	
&+&60b^3c\mu^3\nu+2b^2c^2\mu^2\nu^2-20bc^3\mu\nu^3+63c^4\nu^4\Bigg), \\
		\frac{dc}{dt}&=& \frac{c}{16 a^2 b^2 c^2}\Bigg(63a^4\lambda^4-60a^3b\lambda^3\mu-20a^3c\lambda^3\nu-6a^2b^2\lambda^2\mu^2+20a^2bc\lambda^2\mu\nu \nonumber \\  
&+&2a^2c^2\lambda^2\nu^2-60ab^3\lambda\mu^3+20ab^2c\lambda\mu^2\nu-20abc^2\lambda\mu\nu^2+60ac^3\lambda\nu^3\nonumber \\ 
&+&63b^4\mu^4-20b^3c\mu^3\nu+2b^2c^2\mu^2\nu^2+60bc^3\mu\nu^3-105c^4\nu^4\Bigg). 
	\end{eqnarray}

	{\bf 1: The geometry of $\mathbb{R}+\mathbb{R}+\mathbb{R}$ with the structure constants $\lambda=0 , \mu=0 , \nu=0 $}
	\vspace{0.5cm}
	
	The metric of these geometries are flat, with zero curvature, and since all the flow equations vanish they do not change under the flow. They are fixed points of the flow.
	\begin{equation}
	K^{a}=0, \qquad K=0. \nonumber
	\end{equation}
\vspace{0.5cm}
	
{\bf 2: The geometry of $SU(2)$: with the structure constants $\lambda=1 , \mu=1 , \nu=1 $ or $\lambda=-1 , \mu=-1 , \nu=-1 $}	
\vspace{0.5cm}

For this geometry the flow equations are
	\begin{eqnarray} 
		\frac{da}{dt}&=&\frac{a}{16 a^2 b^2 c^2}\Big(-105a^4+60a^3b+60a^3c+2a^2b^2-20a^2bc+2a^2c^2-20ab^3\nonumber \\ 
&+&20ab^2c+20abc^2-20ac^3+63b^4 -60b^3c-6b^2c^2-60bc^3+63c^4 \Big), \nonumber \\
		\frac{db}{dt}&=& \frac{b}{16 a^2 b^2 c^2}\Big(63a^4-20a^3b-60a^3c+2a^2b^2+20a^2bc-6a^2c^2+60ab^3\nonumber \\ 
&-&20ab^2c+20abc^2-60ac^3-105b^4 +60b^3c+2b^2c^2-20bc^3+63c^4\Big), \nonumber \\
		\frac{dc}{dt}&=& \frac{c}{16 a^2 b^2 c^2}\Big(63a^4-60a^3b-20a^3c-6a^2b^2+20a^2bc+2a^2c^2-60ab^3\nonumber \\ 
&+&20ab^2c-20abc^2+60ac^3+63b^4 -20b^3c+2b^2c^2+60bc^3-105c^4\Big). 
	\end{eqnarray}
	The first observation is that these equations are all symmetric in $a,b$, and $c$: i.e. the interchange of them would not change the equations. In the case of $a(t)= b(t)=c(t)$, the flow equations become rather simple and analytically solvable:
	\begin{equation}
		\frac{da}{dt}=-\frac{1}{16 a} \qquad a^{2}(t)=-\frac{1}{8}t+a_{0}^{2}.
	\end{equation}
	Thus if initially, one has  $a_{0}= b_{0}=c_{0}$, or at some point of the flow $t=\tau$, if $a_\tau= b_\tau=c_\tau$, namely once the metric reaches the metric of the round sphere, the geometry shrinks to a point in a finite time. On the other hand, when $a(t)\neq b(t)\neq c(t)$, it is not possible to solve the equations analytically, so to make some estimates about their behaviors,  let us investigate the differences between the flow equations :
	\begin{eqnarray}
		&\frac{d(a-b)}{dt}&=\frac{(a-b)}{16 a^2 b^2 c^2}\Big(-105a^4-108a^3b+60a^3c-86a^2b^2+100a^2bc+2a^2c^2\nonumber \\ 
                      & -108ab^3&+100ab^2c+28abc^2-20ac^3 -105b^4 +60b^3c+2b^2c^2-20bc^3+63c^4\Big), \nonumber \\ 
		&\frac{d(b-c)}{dt}&=\frac{(b-c)}{16 a^2 b^2 c^2}\Big(63a^4-20a^3b-20a^3c+2a^2b^2+28a^2bc+2a^2c^2+60ab^3  \nonumber \\
		&+100ab^2c&+100abc^2+60ac^3-105b^4-108b^3c-86b^2c^2-108bc^3-105c^4\Big), \nonumber \\ 
		&\frac{d(a-c)}{dt}&=\frac{(a-c)}{16 a^2 b^2 c^2}\Big(-105a^4+60a^3b-108a^3c+2a^2b^2+100a^2bc-86a^2c^2 \nonumber \\
		&-20ab^3&+28ab^2c+100abc^2-108ac^3+63b^4-20b^3c+2b^2c^2+60bc^3-105c^4\Big). 
	\end{eqnarray}
	
	From the symmetry of the flow equations, we can assume without loss of generality that initially $a_{0}\ge b_{0}\ge c_{0}$. If at $t=\tau$, $a_{\tau}=  b_{\tau}$, and similarly for some $t=\tau^{'}$, $b_{\tau^{'}}=  c_{\tau^{'}}$, then one gets
	
	\begin{eqnarray}
		\left.\frac{d(a-b)}{dt}\right \vert_{\tau} &=&\frac{(a-b)}{16 a^2 b^2 c^2}\Big(-105a^4-108a^3 b+60a^3 c-86a^2 b^2+100a^2 bc \nonumber \\
&+&2a^2c^2-108ab^3+100ab^2 c+28abc^2-20ac^3 -105b^4 \nonumber \\
&+&60b^3 c+2b^2 c^2-20bc^3+63c^4\Big ), \nonumber \\
		\left.\frac{d(b-c)}{d t}\right \vert_{\tau^{'}}&=&\frac{(b-c)}{16 a^2 b^2 c^2}\Big(63a^4-20a^3b-20a^3c+2a^2b^2+28a^2b c \nonumber \\ 
&+&2a^2c^2+60ab^3+100ab^2c+100a b c^2+60ac^3 -105b^4 \nonumber \\ 
&-&108b^3c-86b^2c^2-108b c^3-105c^4\Big). 
	\end{eqnarray}
The first equation says that during the flow, the flow between the functions $a$ and $b$ stops at $t=\tau$; and the second equation says that the flow between the functions $b$ and $c$ stops at $t=\tau^{'}$. Therefore if, initially, we have $a_{0}\ge  b_{0}\ge c_{0}$, then for all $t\ge 0$, one has $a\ge  b\ge c$ until they reach exactly the same values. Now, assuming $a\ge  b\ge c$ for all $t$, we want to show that $c$ is a non-decreasing function to the point where all three dimensions become equal. To deal with the negative terms, we write $-60a^3b-60ab^3=30(a-b)^2(a^2-b^2)-30a^4-30b^4$ and $-6a^2b^2=3(a^2-b^2)^2-3a^4-3b^4$. Then one gets 
	\begin{eqnarray}
		\frac{dc}{dt}&=& \frac{c}{16 a^2 b^2 c^2}\Big(10a^4+20a^2bc+2a^2c^2+20ab^2c+60ac^3+2b^2c^2+60bc^3 \nonumber \\
 &+& 10b^4 +3(a^2-b^2)^2+30(a-b)^2(a^2-b^2)-20abc^2-105c^4\Big).
	\end{eqnarray}
	The last two negative terms are smaller than the other terms since $a\ge  b\ge c$. Hence, $c\ge c_{0}$ and so $a\ge  b\ge c\ge c_{0}$. We can put a lower bound for $c$ using the leading term:
	\begin{equation}
		\frac{dc}{dt}= \frac{10a^4c}{16 a^2 b^2 c^2}, \qquad c\frac{dc}{dt}=\frac{5a^2}{8  b^2 } \ge \frac{5}{8,} \qquad c^2(t) \ge \frac{5}{4}t+c_{0}^2.
	\end{equation}
Similarly, we can show that $a$ is non increasing so that we have $a_{0}\ge a\ge  b\ge c\ge c_{0}$. Now we can use this condition for the flow of $a-c$: We see that
	\begin{equation}
	(a_{0}-c_{0})  e^{-\frac{137}{16a_{0}^{2}}t}	\ge (a-c) \ge (a_{0}-c_{0}) e^{-\frac{137}{16c_{0}^{2}}t}.
	\end{equation}
	Therefore $a-c$ vanishes exponentially and similar results follow for others. Hence we conclude that the metric approaches the metric of the round sphere and then the geometry shrinks to a point.
	
	Now let us investigate the behavior of the scalar curvature:  
	\begin{equation}
		R=\frac{1}{2abc}\left(-a^2+2ab+2ac-b^2+2bc-c^2\right).
	\end{equation}
	Once the geometry reaches the round sphere, the curvature scalar becomes $R=\frac{3}{2a}$ and shrinking of the sphere continues; and the scalar curvature becomes singular; it goes to positive infinity.  One can of course avoid this singularity by normalizing the flow. We will not do that here, but for the sake of understanding what goes on in the normalized (fixed volume) flow, let us just note this case as an example.  If we had normalized the flow, then one of  the equations would look like 
\begin{eqnarray}
    	\frac{da}{dt}&=&\frac{a}{6a^2b^2c^2}\Big (-42a^4+25a^3b+25a^3c+a^2b^2-10a^2bc+a^2c^2-5ab^3 \nonumber \\ 
		&+&5ab^2c+5abc^2-5ac^3+21b^4-20b^3c-2b^2c^2-20bc^3+21c^4\Big),
\end{eqnarray}	
and we have similar equations for $b$ and $c$ by interchanging with $a$. In this case, the geometry would reach the round sphere in a finite time and stay there. The scalar curvature then would be $R=3/2a$. Figure  \ref{Fig.SU(2)} shows the behavior of the geometry under the flow (we have done the numerical integrations using a \textsc{Matlab} code).   	
\begin{figure}[h]
    \centering
    \includegraphics[width=120mm,scale=0.5]{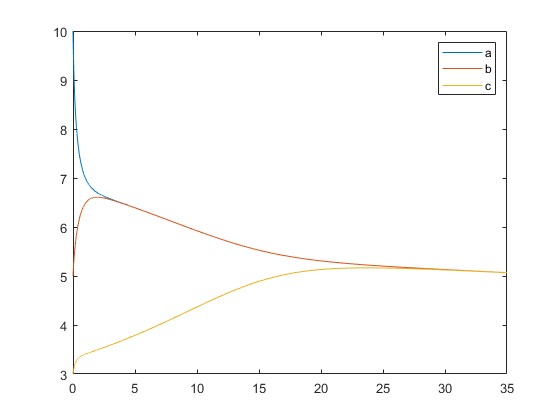}
    \caption{The behavior of the $SU(2) =S^3$ geometry under the flow. The vertical axis represents the dimensions $a,b,c$. The horizontal axis represents the parameter of the flow $t$. Three dimensions approach each other and then vanish in the unnormalized flow, while they attain the same fixed value for the volume normalized flow.}
    \label{Fig.SU(2)}
\end{figure}

When three dimensions become equal we can look at what happens to the trace $K$ which appears as the integrand of the entropy:
\begin{eqnarray}
K &=&\dfrac{1}{32a^2 b^2 c^2}\Big (-21  a^4+20 a^3  b+20  a^3  c+2a^2  b^2-20 a^2  b  c \nonumber \\
&+&2a^2  c^2+20a b^3-20a b^2  c-20 a  b  c^2+20 a  c^3-21 b^4\nonumber \\
&+&20  b^3  c+2 b^2c^2+20  b  c^3-21 c^4\Big )=\dfrac{3}{32a^2 }.
\end{eqnarray}
$K$ first attains a constant value, then as $a$ vanishes, it goes to positive infinity. 
\vspace{0.5cm}
	
For the rest of the Bianchi class geometries we state here the results and show the computations in Appendix B. The geometries of $\widetilde{SL}(2,R)$ and $Nil$ generate a pancake degeneracy. The Geometry of $\widetilde{Isom}(R^{2})$ generates a flat geometry. The geometry of $Sol$ generates a cigar degeneracy. 

\section{Solitons of the $K$-flow}

Fixed points of the flow equations are called solitons. The first kind of solitons of the $K$-flow is the Ricci flat metrics for which the Ricci tensor vanishes. In three dimensions these are just flat spaces.  Under a constant scaling defined as $\tilde{g}_{ij}=\Lambda {g}_{ij}$, the $K$ tensor picks up a $\Lambda^{-1}$ term so scalings are not symmetries of the flow as we have noted before. The second kind of symmetry we may discuss is diffeomorphism.  Letting $\phi_{t}$ be a one-parameter family of diffeomorphisms on the manifold ${\mathcal{M}}$ generated by a vector field $X$ such that $\phi_{0}=id_{\mathcal{M}}$ and $\partial_{t}\phi_{t}=-X(t)$ and a family of metrics 
\begin{equation}
\phi_{t}^{*}g(t)=\bar{g}(t), \qquad g(0)=g_{0}=\bar{g}(0).
\end{equation}
If $\bar{g}(t)$ is the solution of the $K$-flow (\ref{IVP}), then $g(t)$ is a solution of the gauge-fixed flow (\ref{gauged}) by (\ref{Diffeo}). In other words  $g(t)$  is also a solution of the flow with $\partial_{t}g_{ij}=-K_{ij}+\nabla_{i}X_{j}+\nabla_{j}X_{i}$ where we used the definition of the Lie derivative of the metric. We said that a flat metric is a solution of the flow since $K_{ij}=0$ which we may call a {\it steady soliton}.  If $g$ is an Einstein metric ($R_{ij}=2\Lambda g_{ij} $), then $K_{ij}=\frac{1}{2}\Lambda^{2} g_{ij} $ and a conformal Killing vector field $\mathcal{L}_{X}g=\frac{1}{2}\Lambda^{2} g$ would give a soliton solution for an Einstein metric.  We set the pair $(g,X)$ as a soliton of the flow; therefore for a soliton $\nabla_{i}X_{j}+\nabla_{j}X_{i}-K_{ij}=0$. Are there any such solutions? To answer this question,  we start with the square of the Lie derivative of the metric \cite{Cotton2}
\begin{equation}
\mid \mathcal{L}_{X}g \mid^{2} = -\Delta (X_{k}X^{k})+2X_{k} \Delta X^{k}-2\nabla_{k}(X^{l}\nabla_{l}X^{k})+2X^{k}\nabla_{k}\nabla_{l}X^{l}+2X^{k}X^{l}R_{kl},	
\end{equation}
and we also have  $2\nabla_{k}X^{k}-K=0$ and $\Delta X_{k}+\nabla_{k}\nabla_{l}X^{l}+R_{kl}X^{l}=0=\Delta X_{k}+\frac{1}{2}\nabla_{k}K+R_{kl}X^{l}$ from the trace and divergence of the soliton equation. Then the square of the Lie derivative of the metric reads
\begin{equation}
\mid \mathcal{L}_{X}g \mid^{2} = -\Delta (X_{k}X^{k})-2\nabla_{k}(X^{l}\nabla_{l}X^{k}).
\end{equation}
When we integrate over a compact manifold without boundary,  we obtain
\begin{equation}
\int_{\mathcal{M}} d\mu \mid \mathcal{L}_{X}g \mid^{2} = 0.
\end{equation}
Apparently, only a Killing vector field solves this equation. Then soliton equation gives a trivial soliton $K_{ij}=0$. In the case of $G$-invariant metrics, the left-invariant vector fields are Killing and are generators of diffeomorphisms, in fact, isometries, and the associated integral curves are geodesics of the manifold. Steady solitons of the flow, flat metrics, flow on the geodesics. 

From previous discussions, Einstein metrics are critical points of the action, and so of the flow. We should emphasize that Einstein metrics as critical points of the quadratic part of the NMG action are found from purely variational methods in \cite{Viac}. As it stands we cannot require a general solution of the flow to vanish at the variations of the metric. Hence we conclude that Einstein metrics are critical points, but not solitons of the flow.

\section{Conclusions}

We have defined a geometric flow exclusive to 3-dimensional Riemannian manifolds. The flow is motivated by both pure differential geometry and physics. Compared to the Ricci flow, the $ K$ flow is a higher-order flow, as the relevant tensor involves the laplacian of the Ricci tensor at the leading order. Despite being a higher-order non-linear PDE, we were able to show that the relevant operator is a fourth-order elliptic operator using the deTuck trick and establishing the short-time existence problem. The fast modes (short wavelength perturbations) are all diffused just as in the case of the Ricci flow. From the physics point of view, the relevant $K$-tensor can be thought to arise from the short distance (ultra-violate regime)  of the new massive gravity for which the Ricci tensor (or the Einstein tensor) becomes an irrelevant operator \cite{nmg1}, or one can directly consider the quadratic theory as the gravity theory as was done in\cite{Deser}. From a mathematical point of view, what we called ``the entropy'' of the flow was studied in \cite{Viac} in the context of characterizing 3-dimensional space forms. We have made use of the results of  \cite{Viac} in section III where we gave a gradient formulation of the flow. 

Ricci flow and various related flows bore much fruit in mathematics as is already well known from the works of Hamilton and Perelman that we alluded to in the introduction. But these remarkable works, and in general works on geometric flows have not yet found their deserved interest in the physics literature. We believe this state of affairs is transient: There are a great deal of potential applications of geometric flows in gravity, quantum gravity, and field theory. The crucial link is the following: Ricci flow, or some nice geometric flow, could arise as a renormalization group flow in a non-linear sigma model. For example when one says that Einstein's theory arises from string theory, one means the following: at low energies, the worldsheet renormalization group flow in string theory gives rise to the Ricci flow (with some additional fields) in the target space. And, because in the world sheet, Weyl symmetry cannot be anomalous, one sets the beta function (that is running of the couplings which is the metric of the target space) to zero. So instead of a flow in string theory, one has the critical points of the flow which are at the lowest order, and setting all the fields except the metric to vanish, to the Ricci flat metrics.  Of course, as we shall mention below, in the context of Euclidean gravity, Ricci flow could still be useful as the natural path connecting the critical (or saddle) points. 

 Let us give some examples from the recent literature of the use of geometric flows. In \cite{Gegenberg} a modified Ricci flow with Maxwell-Chern-Simons terms and additional fields, inspired by string theory was studied within the context of the $3D$ uniformization theorem.  In \cite{Samuel}, to find inequalities in General Relativity regarding the evolution of the area of a surface and the enclosed Hawking mass, Ricci flow was employed.  In \cite{Cotton1, Cotton2} Cotton flow was studied to understand the flow of homogenous geometries. 100\% Cotton flow arises as the high energy limit of the Horava gravity \cite{Horava} and the topologically massive gravity \cite{TMG} which were studied in 
\cite{Bakas1, Bakas2}  and \cite{Maloney} respectively. In \cite{Headrick}, Ricci flow was used to study the transitions between saddle points of the four-dimensional Euclidean gravity. Ricci flow appears as a gradient flow between the possible saddle points. Of course, there could be many such saddle points, but once a boundary topology is introduced, one can handle a few saddle points. In the work, the boundary was taken to be $S^1\times S^2$, and with this boundary, there are 3 saddle points. In \cite{Woolgar} a nice compilation of possible applications of the Ricci flow in physics was given. Recently, in \cite{Frenkel} a topological non-relativistic quantum gravity was defined in which the Ricci flow equation appears as the localization formula in the path integral that defines the theory. 

In this work, we laid the basic structure of the $K$-flow; and applied it to the homogenous geometries that are potentially relevant to the spatial geometry of the universe at large scales. Applications of the $K$-flow  for non-homogenous geometries, such as black holes are also of interest. Various such solutions as the critical points of this flow (i.e. those that satisfy $K_{i j}=0$) have been found in \cite{Barnich} and elaborated in \cite{Ercan2}. These are black holes with deformed horizons (also called black flowers), static black holes, rotating black holes, and dynamical black flowers that emit or absorb gravitons.  This set of critical solutions constitutes a potentially very useful application to the $K$-flow. For example, the Euclidean version of these solutions will be the saddle points of the flow, and one can study the transition between these solutions as was done in \cite{Headrick} for the Ricci flow. Moreover, $K$-flow can be considered as the high energy limit of the ``new massive gravity'' inspired flow with $ \partial_t g_{\i j} = \alpha K_{\i j} + \beta G_{i j} + \Lambda g_{\i j}$ of which the critical points are metrics that solve the new massive gravity \cite{nmg1,nmg2} that has been studied as a dynamical theory of gravity  with a massive spin-2 graviton in three dimensions (as opposed to the locally trivial Einstein's gravity in three dimensions). Many solutions of this extended theory are known which can be studied as saddle points under the flow.

\section{Appendices}

In these appendices, we give some details of the computations that we used in the text and also compute the properties of the $K$-tensor under conformal transformations. 

\subsection{ Computation of the $K$-tensor in Orthonormal Frame}

We now compute the Ricci tensor and the curvature scalar  following (\ref{EFmilnor}) 
\begin{align}
[E_{1},E_{2}]=\nu \sqrt{\dfrac{c}{ab}}E_{3} , \qquad [E_{2},E_{3}]=\lambda\sqrt{\dfrac{a}{bc}} E_{1}  , \qquad [E_{3},E_{1}]=\mu \sqrt{\dfrac{b}{ac}} E_{2} \quad .
\end{align}

\begin{align}
\nabla_{E_{1}}E_{2}=\dfrac{1}{2}( D _{12}~^{3} -D _{23}~^{1}+D _{31}~^{2} ) E_{3} =\dfrac{1}{2}\left( \nu \sqrt{\dfrac{c}{ab}}  -  \lambda \sqrt{\dfrac{a}{bc}} + \mu \sqrt{\dfrac{b}{ac}} \right)  E_{3} \quad .
\end{align}
Therefore one obtains 
\begin{align}
\nabla_{F_{1}}F_{2}&=\dfrac{1}{2}\left( \dfrac{-\lambda a+\mu b +\nu c}{c}\right)  F_{3}\quad , \qquad
\nabla_{F_{2}}F_{1}=\dfrac{1}{2}\left( \dfrac{-\lambda a+\mu b -\nu c}{c}\right)  F_{3} \quad ,  \nonumber \\
\nabla_{F_{1}}F_{3}&=\dfrac{1}{2}\left( \dfrac{\lambda a-\mu b -\nu c}{b}\right)  F_{2}\quad , \qquad
\nabla_{F_{3}}F_{1}=\dfrac{1}{2}\left( \dfrac{\lambda a+\mu b -\nu c}{b}\right)  F_{2}\quad ,  \nonumber \\
\nabla_{F_{2}}F_{3}&=\dfrac{1}{2}\left( \dfrac{\lambda a-\mu b +\nu c}{a}\right)  F_{1}\quad , \qquad
\nabla_{F_{3}}F_{2}=\dfrac{1}{2}\left( \dfrac{-\lambda a-\mu b +\nu c}{a}\right)  F_{1}\quad .
 \end{align}
One can find the sectional curvatures and then the components of the Ricci tensor by the use of these equations. The usefulness of this frame is that the Ricci tensor is diagonal which is very handy when studying the flows. The non-zero diagonal components of the Ricci tensor in the orthonormal frame can be computed as 
\begin{align}
R_{11}&=\dfrac{1}{2abc}\left(\lambda^{2}a^{2}-\mu^{2}b^{2}-\nu^{2}c^{2}+2\mu\nu bc\right) , \qquad
R_{22}=\dfrac{1}{2abc}\left(-\lambda^{2}a^{2}+\mu^{2}b^{2}-\nu^{2}c^{2}+2\lambda\nu ac\right) ,  \nonumber \\
R_{33}&=\dfrac{1}{2abc}\left(-\lambda^{2}a^{2}-\mu^{2}b^{2}+\nu^{2}c^{2}+2\lambda\mu ab\right)  . 
\end{align}
The non-zero diagonal components of the Ricci tensor in the orthogonal frame by the relation between two coordinates can be computed as
\begin{align}
R_{11}&=\dfrac{a}{2abc}\left(\lambda^{2}a^{2}-\mu^{2}b^{2}-\nu^{2}c^{2}+2\mu\nu bc\right) , \qquad
R_{22}=\dfrac{b}{2abc}\left(-\lambda^{2}a^{2}+\mu^{2}b^{2}-\nu^{2}c^{2}+2\lambda\nu ac\right) ,  \nonumber \\
R_{33}&=\dfrac{c}{2abc}\left(-\lambda^{2}a^{2}-\mu^{2}b^{2}+\nu^{2}c^{2}+2\lambda\mu ab\right) . 
\end{align}
The scalar curvature  then will be :
\begin{align}
R=\dfrac{1}{2abc}\left(-\lambda^{2}a^{2}-\mu^{2}b^{2}-\nu^{2}c^{2}+2\lambda\mu ab+2\lambda\nu ac+2\mu\nu bc\right) \quad .
\end{align}
We now find the connections and the curvature form:
\begin{align}
de^{1}&=d(\sqrt{a}\omega^{1})=\sqrt{a}d\omega^{1}=\sqrt{a} \lambda\omega^{2}\wedge\omega^{3}=\sqrt{\dfrac{a}{bc}} \lambda e^{2}\wedge e^{3} \quad ,  \nonumber \\
de^{2}&=d(\sqrt{b}\omega^{2})=\sqrt{b}d\omega^{2}=\sqrt{b} \mu\omega^{3}\wedge\omega^{1}=\sqrt{\dfrac{b}{ac}} \mu e^{3}\wedge e^{1} \quad ,  \nonumber \\
de^{3}&=d(\sqrt{c}\omega^{3})=\sqrt{c}d\omega^{3}=\sqrt{c} \nu\omega^{1}\wedge\omega^{2}=\sqrt{\dfrac{c}{ab}} \nu e^{1}\wedge e^{2} \quad .
\end{align}
In the orthonormal frame $g_{ab}=\delta_{ab}$, hence the metric compatibility is $\nabla_{c}\delta_{ab}=0$ which yields  $\omega^{b}~_{a}=-\omega^{a}~_{b}$. Then one has 
\begin{align}
\omega^{1}~_{2}&=\dfrac{1}{2\sqrt{abc}}\left(\lambda a+\mu b- \nu c  \right)e^{3}, \qquad
\omega^{1}~_{3}=\dfrac{1}{2\sqrt{abc}}\left(-\lambda a+\mu b- \nu c  \right)e^{2}, \nonumber \\
\omega^{2}~_{3}&=\dfrac{1}{2\sqrt{abc}}\left(-\lambda a+\mu b+\nu c  \right)e^{1} .
\end{align}
The components of the curvature 2-form can be computed by the relation $R^{a}~_{b}=d\omega^{a}~_{b}+\omega^{a}~_{c}\wedge \omega^{c}~_{b}$ .  For example, 
\begin{align}
R^{1}~_{2}=d\omega^{1}~_{2}+\omega^{1}~_{1}\wedge \omega^{1}~_{2}+\omega^{1}~_{2}\wedge \omega^{2}~_{2}+\omega^{1}~_{3}\wedge \omega^{3}~_{2}=d\omega^{1}~_{2}+\omega^{1}~_{3}\wedge \omega^{3}~_{2} \quad ,
\end{align}
 where
\begin{align}
 d\omega^{1}~_{2}&=\dfrac{1}{2\sqrt{abc}}\left(\lambda a+\mu b- \nu c \right)de^{3}=\dfrac{1}{2\sqrt{abc}}\left(\lambda a+\mu b- \nu c  \right)\sqrt{\dfrac{c}{ab}} \nu e^{1}\wedge e^{2} \nonumber \\
&=\dfrac{\nu}{2ab}\left(\lambda a+\mu b- \nu c \right)  e^{1}\wedge e^{2},
\end{align}
\begin{align}
\omega^{1}~_{3}\wedge \omega^{3}~_{2}&=\Big(\dfrac{1}{2\sqrt{abc}}\left(-\lambda a+\mu b- \nu c  \right)e^{2}\Big)\wedge \Big(-\dfrac{1}{2\sqrt{abc}}\left(-\lambda a+\mu b+\nu c  \right)e^{1}\Big), \nonumber \\
&=\dfrac{1}{4abc}\Big(\lambda^{2} a^{2}+\mu^{2} b^{2}-\nu^{2} c^{2}-2\lambda \mu ab  \Big)e^{1}\wedge e^{2}.
\end{align}
These computations lead to 
\begin{align}
R^{1}~_{2}&=\dfrac{1}{4abc}\Big(\lambda^{2} a^{2}+\mu^{2} b^{2}-3\nu^{2} c^{2}-2\lambda \mu ab+2\lambda \nu ac +2\mu \nu bc  \Big)e^{1}\wedge e^{2}, \nonumber \\
R^{1}~_{3}&=\dfrac{1}{4abc}\Big(\lambda^{2} a^{2}-3\mu^{2} b^{2}+\nu^{2} c^{2}+2\lambda \mu ab-2\lambda \nu ac +2\mu \nu bc  \Big)e^{1}\wedge e^{3},  \nonumber \\
R^{2}~_{3}&=\dfrac{1}{4abc}\Big(-3\lambda^{2} a^{2}+\mu^{2} b^{2}+\nu^{2} c^{2}+2\lambda \mu ab+2\lambda \nu ac -2\mu \nu bc \Big)e^{2}\wedge e^{3}.
\end{align}
Ricci 1-forms can be computed by the relation $Ric_{a}=\imath_{b}R^{b}~_{a}$. For example 
\begin{align}
 Ric_{1}&=\imath_{2}R^{2}~_{1}+\imath_{3}R^{3}~_{1} \nonumber \\
&=\dfrac{1}{4abc}\Big(\lambda^{2} a^{2}+\mu^{2} b^{2}-3\nu^{2} c^{2}-2\lambda \mu ab+2\lambda \nu ac +2\mu \nu bc \Big)\imath_{2}(e^{2}\wedge e^{1}) \nonumber \\
& + \dfrac{1}{4abc}\Big(\lambda^{2} a^{2}-3\mu^{2} b^{2}+\nu^{2} c^{2}+2\lambda \mu ab-2\lambda \nu ac +2\mu \nu bc \Big)\imath_{3}(e^{3}\wedge e^{1}), \nonumber \\
\end{align}
where the interior product acts as
\begin{align}
\imath_{2}(e^{2}\wedge e^{1})&=(\imath_{2}e^{2}) \wedge e^{1}+(-1)^{1}e^{2}\wedge (\imath_{2}e^{1})=e^{1}, \nonumber \\
\imath_{3}(e^{3}\wedge e^{1})&=(\imath_{3}e^{3}) \wedge e^{1}+(-1)^{1}e^{3}\wedge (\imath_{3}e^{1})=e^{1}.
\end{align}
Hence we find Ricci 1-forms as 
\begin{align}
Ric_{1}&=\dfrac{1}{2abc}\left(\lambda^{2} a^{2}-\mu^{2} b^{2}-\nu^{2} c^{2} +2\mu \nu bc \right)e^{1}, \nonumber \\
Ric_{2}&=\dfrac{1}{2abc}\left(-\lambda^{2} a^{2}+\mu^{2} b^{2}-\nu^{2} c^{2} +2\lambda \nu ac\right)e^{2}, \nonumber \\
Ric_{3}&=\dfrac{1}{2abc}\left(-\lambda^{2} a^{2}-\mu^{2} b^{2}+\nu^{2} c^{2} +2\lambda  \mu ab\right)e^{3}.
\end{align}
The curvature scalar can be computed from the equation $R=\imath_{a}(Ric)^{a}$ as
\begin{align}
R=\dfrac{1}{2abc}\left(-\lambda^{2} a^{2}-\mu^{2} b^{2}-\nu^{2} c^{2} +2\lambda  \mu ab+2\lambda\nu ac+ 2\mu\nu bc\right).
\end{align}
Schouten one forms easily computed from the earlier results as
$S^{a}=(Ric)^{a}-\dfrac{1}{4}R e^{a}$. Namely
\begin{align}
S^{1}&= \dfrac{1}{2abc} \left(\frac{5\lambda^{2} a^{2}}{4}-\frac{\lambda a\mu b}{2}-\frac{\lambda a\nu c}{2}-\frac{3\mu^{2} b^2}{4}+\frac{3\mu b\nu c}{2}-\frac{3\nu^{2} c^2}{4}\right), \nonumber \\
S^{2}&=\dfrac{1}{2abc} \left(-\frac{3\lambda^{2} a^2}{4}-\frac{\lambda a\mu b}{2}+\frac{3\lambda a \nu c}{2}+\frac{5\mu^{2} b^2}{4}-\frac{\mu b\nu c}{2}-\frac{3\nu^{2} c^2}{4}\right), \nonumber \\
S^{3}&=\dfrac{1}{2abc} \left(-\frac{3\lambda^{2} a^2}{4}+\frac{3\lambda a\mu b\nu c}{2}-\frac{\lambda a\nu c}{2}-\frac{3\mu^{2} b^2}{4}-\frac{\mu b\nu c}{2}+\frac{5\nu^{2} c^2}{4}\right). 
\end{align}
The Cotton 2-form is in the following form: 
\begin{equation}
 C^{a}=DS^{a}=d\Big ( (Ric)^{a}-\dfrac{1}{4}R e^{a}\Big)+\omega^{a}~_{b}\wedge \Big( (Ric)^{b}-\dfrac{1}{4}R e^{b}\Big) . 
\end{equation}
Hence 
\begin{align}
C^{1}&=\dfrac{1	}{2(abc)^{3/2}}\Big(-\lambda^{2}a^{2}(-2\lambda a+\mu b+\nu c)-(\mu b+\nu c)( \mu b - \nu c )^{2}\Big) e^{2}\wedge e^{3}, \nonumber \\
C^{2}&=\dfrac{1	}{2(abc)^{3/2}}\Big(-\mu^{2}b^{2}( \lambda a-2\mu b +\nu c)-(\lambda  a+\nu c)( \lambda a - \nu c )^{2}\Big) e^{3}\wedge e^{1}, \nonumber \\
C^{3}&=\dfrac{1	}{2(abc)^{3/2}}\Big( -\nu^{2}c^{2}( \lambda a+\mu b -2\nu c)-(\lambda a +\mu b)( \lambda a - \mu b )^{2}\Big) e^{1}\wedge e^{2}.
\end{align}
We find the Cotton one forms as the Hodge dual of the Cotton 2-forms $\ast C^{a}$ and arrive at
\begin{align}
\ast C^{1}&=\dfrac{1}{2(abc)^{3/2}}\Big(-\lambda^{2}a^{2}(-2\lambda a+\mu b+\nu c)-(\mu b+\nu c)( \mu b - \nu c )^{2}\Big)e^{1}, \nonumber \\
\ast C^{2}&=\dfrac{1}{2(abc)^{3/2}}\Big(-\mu^{2}b^{2}( \lambda a-2\mu b +\nu c)-(\lambda  a+\nu c)( \lambda a - \nu c )^{2}\Big) e^{2}, \nonumber \\
\ast C^{3}&=\dfrac{1}{2(abc)^{3/2}}\Big(-\nu^{2}c^{2}( \lambda a+\mu b -2\nu c)-(\lambda a +\mu b)( \lambda a - \mu b )^{2}\Big) e^{3}. 
\end{align}

Finally we can compute the one-forms $J^a=\frac{1}{2}\epsilon_{\phantom{a}bc}^{a}\star (S^{b}\wedge S^{c})$: 
\begin{eqnarray}
J^{1}&=& \dfrac{1}{4a^2 b^2 c^2}	\Bigg(\frac{9\lambda^{4} a^4}{16}-\frac{3\lambda^{3} a^3 \mu b}{4}-\frac{3\lambda^{3} a^3 \nu c}{4}-\frac{9\lambda^{2} a^2 \mu^{2} b^2}{8}+\frac{13\lambda^{2} a^2 \mu b\nu c}{4} \nonumber \\ 
&-&\frac{9\lambda^{2} a^2 \nu^{2} c^2}{8} +\frac{9\lambda a\mu^{3} b^3}{4}-\frac{9\lambda a \mu^{2} b^2 \nu c}{4} -\frac{9\lambda a\mu b\nu^{2} c^2}{4}+\frac{9\lambda a \nu^{3} c^3}{4}-\frac{15\mu^{4} b^4}{16}\nonumber \\
&-&\frac{\mu^{3}  b^3  \nu c}{4}+\frac{19\mu^{2} b^2 \nu^{2} c^2}{8}-\frac{\mu b\nu^{3} c^3}{4}-\frac{15\nu^{4} c^4}{16}\Bigg)e^{1}, \nonumber \\
J^{2}&=&\dfrac{1}{4a^2 b^2 c^2}\Bigg (-\frac{15\lambda^{4} a^4}{16}+\frac{9\lambda^{3} a^3 \mu b}{4}-\frac{\lambda^{3} a^3 \nu c}{4}-\frac{9\lambda^{2} a^2 \mu^{2} b^2}{8}-\frac{9\lambda^{2} a^2 \mu b \nu c}{4} \nonumber \\ 
&+&\frac{19\lambda^{2} a^2 \nu^{2} c^2}{8} -\frac{3\lambda a \mu^{3} b^3}{4}+\frac{13\lambda a \mu^{2} b^2 \nu c}{4}-\frac{9\lambda a \mu b \nu^{2} c^2}{4}-\frac{ \lambda a \nu^{3} c^3}{4}+\frac{9\mu^{4} b^4}{16}\nonumber \\
&-&\frac{3 \mu^{3} b^3 \nu c}{4}-\frac{9 \mu^{2} b^2 \nu^{2} c^2}{8}+\frac{9 \mu b \nu^{3} c^3}{4}-\frac{15 \nu^{4} c^4}{16}\Bigg)e^{2}, \nonumber \\
J^{3}&=&\dfrac{1}{4a^2 b^2 c^2}\Bigg (-\frac{15\lambda^{4} a^4}{16}-\frac{\lambda^{3} a^3 \mu b}{4}+\frac{9\lambda^{3} a^3 \nu c}{4}+\frac{19\lambda^{2} a^2 \mu^{2} b^2}{8}-\frac{9\lambda^{2} a^2 \mu b \nu c}{4} \nonumber \\ 
&-&\frac{9\lambda^{2} a^2 \nu^{2} c^2}{8} -\frac{\lambda a \mu^{3} b^3}{4}-\frac{9\lambda a \mu^{2} b^2 \nu c}{4}+\frac{13\lambda a \mu b \nu^{2} c^2}{4}-\frac{3 \lambda a \nu^{3} c^3}{4}-\frac{15\mu^{4} b^4}{16} \nonumber \\ 
&+&\frac{9 \mu^{3} b^3 \nu c}{4}-\frac{9 \mu^{2} b^2 \nu^{2} c^2}{8}-\frac{3 \mu b \nu^{3} c^3}{4}+\frac{9 \nu^{4} c^4}{16}\Bigg)e^{3}. 
\end{eqnarray}
The $H^a=\star D\star C^{a}$ one forms read:
\begin{eqnarray}
H^{1}&=& \dfrac{1}{4a^2 b^2 c^2}	\Bigg(6\lambda^4 a^4-3\lambda^3 a^3 \mu b-3\lambda^3 a^3\nu c+\lambda^2 a^2 \mu^2 b^2-2\lambda^2 a^2 \mu b \nu c + \lambda^2 a^2 \nu^2 c^2 \nonumber \\
&-&\lambda a \mu^3 b^3+\lambda a \mu^2 b^2 \nu c+\lambda a \mu b \nu^2 c^2-\lambda a \nu^3 c^3-3\mu^4 b^4+4\mu^3 b^3 \nu c-2\mu^2  b^2 \nu^2 c^2\nonumber \\
&+&4\mu b \nu^3 c^3-3\nu^4 c^4\Bigg)e^{1}, \nonumber \\
H^{2}&=& \dfrac{1}{4a^2 b^2 c^2}	\Bigg(-3\lambda^4 a^4-\lambda^3 a^3 \mu b+4\lambda^3 a^3\nu c+\lambda^2 a^2 \mu^2 b^2+\lambda^2 a^2 \mu b \nu c -2 \lambda^2 a^2 \nu^2 c^2 \nonumber \\
&-&3\lambda a \mu^3 b^3-2\lambda a \mu^2 b^2 \nu c+\lambda a \mu b \nu^2 c^2+4\lambda a \nu^3 c^3+6\mu^4 b^4-3\mu^3 b^3 \nu c+\mu^2  b^2 \nu^2 c^2 \nonumber \\
&-&\mu b \nu^3 c^3-3\nu^4 c^4\Bigg)e^{2}, \nonumber \\
H^{3}&=& \dfrac{1}{4a^2 b^2 c^2}	\Bigg(-3\lambda^4 a^4+4\lambda^3 a^3 \mu b-\lambda^3 a^3\nu c-2\lambda^2 a^2 \mu^2 b^2+\lambda^2 a^2 \mu b \nu c + \lambda^2 a^2 \nu^2 c^2 \nonumber \\
&+&4\lambda a \mu^3 b^3+\lambda a \mu^2 b^2 \nu c-2\lambda a \mu b \nu^2 c^2-3\lambda a \nu^3 c^3-3\mu^4 b^4-\mu^3 b^3 \nu c+\mu^2  b^2 \nu^2 c^2\nonumber \\
&-&3\mu b \nu^3 c^3+6\nu^4 c^4\Bigg)e^{3}. 
\end{eqnarray}

$H$ is indeed traceless. We finally find the $K$ one forms to be

	\begin{align}
K^{1}&= -\frac{1}{32 a^2 b^2 c^2}\Bigg(-105a^4\lambda^4+60a^3b\lambda^3\mu+60a^3c\lambda^3\nu+2a^2b^2\lambda^2\mu^2-20a^2bc\lambda^2\mu\nu \nonumber \\ 
 &+2a^2c^2\lambda^2\nu^2 -20ab^3\lambda\mu^3+20ab^2c\lambda\mu^2\nu+20abc^2\lambda\mu\nu^2-20ac^3\lambda\nu^3+63b^4\mu^4 \nonumber \\ 
&-60b^3c\mu^3\nu-6b^2c^2\mu^2\nu^2-60bc^3\mu\nu^3+63c^4\nu^4\Bigg) e^{1}, \nonumber \\
K^{2}&=- \frac{1}{32 a^2 b^2 c^2} \Bigg(63a^4 \lambda^4-20a^3 b \lambda^3 \mu -60a^3 c \lambda^3 \nu + 2a^2 b^2 \lambda^2 \mu^2 +20a^2 b c \lambda^2 \mu \nu \nonumber \\ 
 &-6a^2 c^2 \lambda^2 \nu^2+60ab^3 \lambda \mu^3 -20ab^2 c \lambda \mu^2 \nu +20abc^2 \lambda \mu \nu^2 -60ac^3 \lambda \nu^3-105b^4 \mu^4 \nonumber \\
&+60b^3 c \mu^3 \nu +2b^2 c^2 \mu^2 \nu^2-20bc^3\mu\nu^3+63c^4\nu^4\Bigg)e^{2}, \nonumber \\
K^{3}&= -\frac{1}{32 a^2 b^2 c^2}\Bigg(63a^4\lambda^4-60a^3b\lambda^3\mu-20a^3c\lambda^3\nu-6a^2b^2\lambda^2\mu^2+20a^2bc\lambda^2\mu\nu \nonumber \\  
&+2a^2c^2\lambda^2\nu^2-60ab^3\lambda\mu^3+20ab^2c\lambda\mu^2\nu-20abc^2\lambda\mu\nu^2+60ac^3\lambda\nu^3\nonumber \\  
&+63b^4\mu^4-20b^3c\mu^3\nu+2b^2c^2\mu^2\nu^2+60bc^3\mu\nu^3-105c^4\nu^4\Bigg) e^{3}. 
	\end{align}
	
The trace of the $K$ tensor turns out to be 
\begin{eqnarray}
K &=&\dfrac{1}{2a^2 b^2 c^2}\Bigg(-\frac{21 \lambda^4 a^4}{16}+\frac{5 \lambda^3 a^3 \mu b}{4}+\frac{5 \lambda^3 a^3 \nu c}{4}+\frac{\lambda^2 a^2 \mu^2 b^2}{8}-\frac{5\lambda^2 a^2 \mu b \nu c}{4} \nonumber \\
&+&\frac{\lambda^2 a^2 \nu^2 c^2}{8}+\frac{5 \lambda a \mu^3 b^3}{4}-\frac{5\lambda a \mu^2 b^2 \nu c}{4}-\frac{5\lambda a \mu b \nu^2 c^2}{4}+\frac{5 \lambda a \nu^3 c^3}{4}-\frac{21 \mu^4 b^4}{16}\nonumber \\
&+&\frac{5 \mu^3 b^3 \nu c}{4}+\frac{\mu^2 b^2 \nu^2 c^2}{8}+\frac{5 \mu b \nu^3 c^3}{4}-\frac{21 \nu^4 c^4}{16}\Bigg).
\end{eqnarray}
Of course we could compute the trace from the Lagrangian $\frac{3}{8}R^{2}\star 1 -R^{a}\wedge \star R_{a} $. 

\subsection{The Bianchi class geometries continued}
	{\bf 3: The geometry of $\widetilde{SL}(2,R)$ with the structure constants $\lambda=-1 , \mu=-1 , \nu=1 $}
	\vspace{0.5cm}
	
For this geometry the flow equations are
	\begin{eqnarray} 
		\frac{da}{dt}&=&\frac{a}{16 a^2 b^2 c^2}\Big (-105a^4+60a^3b-60a^3c+2a^2b^2+20a^2bc+2a^2c^2-20ab^3\nonumber \\
		 &-&20ab^2c+20abc^2+20ac^3+63b^4+60b^3c-6b^2c^2+60bc^3+63c^4\Big),  \nonumber \\
		\frac{db}{dt}&=& \frac{b}{16 a^2 b^2 c^2}\Big(63a^4-20a^3b+60a^3c+2a^2b^2-20a^2bc-6a^2c^2+60ab^3\nonumber \\
		 &+&20ab^2c+20abc^2+60ac^3-105b^4-60b^3c+2b^2c^2+20bc^3+63c^4\Big), \nonumber\\
		\frac{dc}{dt}&=& \frac{c}{16 a^2 b^2 c^2}\Big (63a^4-60a^3b+20a^3c-6a^2b^2-20a^2bc+2a^2c^2-60ab^3\nonumber \\ 
		&-&20ab^2c-20abc^2-60ac^3+63b^4+20b^3c+2b^2c^2-60bc^3-105c^4\Big). 
	\end{eqnarray}
These equations are symmetric in $a$ and $b$, so we may set $a_{0}\ge b_{0}$. If at any $t=\tau$,   $a_{\tau}= b_{\tau}$
	\begin{eqnarray}
		\left.\frac{d(a-b)}{d t}\right\vert_{t=\tau}&=& \frac{a-b}{16 a^2 b^2 c^2}\Big(-105a^4-108a^3b-60a^3c-86a^2b^2-100a^2b c \nonumber \\   
		&+&2a^2c^2-108ab^3-100ab^2c+28a b c^2+20ac^3 \nonumber \\  
		&-&105b^4-60b^3c+2b^2c^2+20b c^3+63c^4\Big). 
	\end{eqnarray}
	Thus if initially $a_{0}\ge b_{0}$, then $a\ge b$ for all $t\ge 0$. From the antisymmetry between $a$ and $c$, we assume that they behave oppositely under the flow. If at some time, all three dimensions are equal, namely if they attain this, then the equation for $a$ would be 
	\begin{equation}
\frac{da}{d t}=\frac{159 a}{16}	\rightarrow a^{2}(t)=a_{0}^{2}+\frac{159}{8} t.
\end{equation}	 
	Therefore we conclude that $a$ and $b$ are increasing whereas $c$ is a decreasing function of $t$ generating a pancake degeneracy under the flow. At a finite time $c$ vanishes, and when that happens, since $a=b$, the flows of all three dimensions stop. The figure  \ref{Fig.SL(2,R)} shows the behavior of the geometry under the flow.	
\begin{figure}[h]
    \centering
    \includegraphics[width=120mm,scale=0.5]{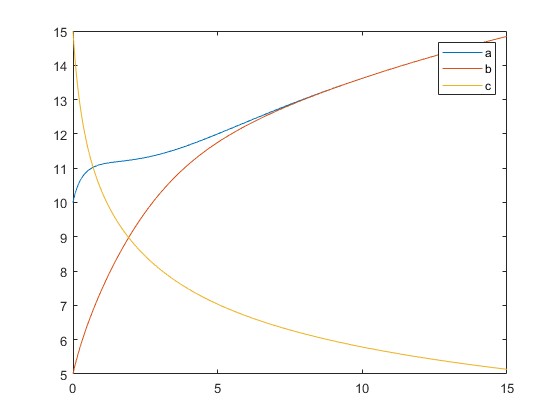}
    \caption{The behavior of the $SL(2,R)$ geometry under the flow. The vertical axis represents the dimensions $a,b,c$. The horizontal axis represents the parameter of the flow $t$. Two dimensions $a,b$ expand, while $c$ vanishes under the flow, yielding a pancake degeneracy.}
    \label{Fig.SL(2,R)}
\end{figure}
	
Next, let us investigate the scalar curvature.  
	\begin{eqnarray}
		R=\frac{1}{2abc}\left(-(a-b)^2-c^2-2c(a+b)\right).
	\end{eqnarray}	
As $a,b$ approach to each other, $R$ attains a negative constant value $(-4ac-c^{2})/2a^{2}c$. Then by the expansion of $a,b$ and the vanishing of $c$, $R$ vanishes. On the other hand, the quadratic curvature invariant
\begin{eqnarray}
K &=&\dfrac{1}{32a^2 b^2 c^2}\Big(-21  a^4+20  a^3 b-20 a^3 c+2 a^2  b^2+20 a^2  b  c +2 a^2  c^2+20  a  b^3 \nonumber\\
&+&20 a b^2 c-20 a  b  c^2-20  a  c^3-21 b^4-20  b^3  c+2 b^2 c^2-20 b c^3-21  c^4\Big)
\end{eqnarray}	
attains a negative constant value $K=(-16a^{2}-40ac-21c^{2})/32a^{4}$ when $a=b$; and when $c$ vanishes, it becomes $K=-1/2a^{2}$ and stays constant.
\vspace {0.5 cm}
 
{\bf 4:  The Geometry of $\widetilde{Isom}(R^{2})$ with the structure constants $\lambda=-1 , \mu=-1 , \nu=0 $}
	\vspace {0.5 cm}
	
	For this geometry, the flow equations are
		\begin{eqnarray} 
		\frac{da}{dt}&=&\frac{a}{16 a^2 b^2 c^2}\Big(-105a^4+60a^3b+2a^2b^2-20ab^3+63b^4\Big),  \nonumber \\
		\frac{db}{dt}&=& \frac{b}{16 a^2 b^2 c^2}\Big(63a^4-20a^3b+2a^2b^2+60ab^3-105b^4\Big), \nonumber \\
		\frac{dc}{dt}&=& \frac{c}{16 a^2 b^2 c^2}\Big(63a^4-60a^3b-6a^2b^2-60ab^3+63b^4\Big).
	\end{eqnarray}
	These equations are symmetric in $a$ and $b$, so let us look at the difference between these two dimensions:
	\begin{equation}
\frac{d(a-b)}{dt}=\frac{(a-b)}{16 a^2 b^2 c^2}\Big(-105a^4-108a^3b-86a^2b^2-108ab^3-105b^4\Big).
	\end{equation}
	If at some time $t=\tau$, $a=b$ then the flow stops. Thus if initially $a_{0}\geq b_{0}$ then for all times one has $a \geq b$, and this implies also that $(a-b)$ is a decreasing function. Moreover if ever $a=b$, then the flow of $c$ stops too.  Using this condition on $a$ and $b$, we can make some estimates about the behavior of the flow equations. After some manipulations, we get
	\begin{eqnarray}
	\frac{da}{dt}&=&\frac{a}{16 a^2 b^2 c^2}\Bigg(-73(a^4-b^4)-2a^2(a^2-b^2)-20b^3(a-b) \nonumber \\ 
	&-&3(a^2-b^2)(a-b)^2 -60ab^3\Bigg) \leq -\frac{15b}{4c^2}, \nonumber \\
	\frac{db}{dt}&=& \frac{b}{16 a^2 b^2 c^2}\Bigg(53(a^4-b^4)+2b^2(a^2-b^2)+60b^3(a-b) \nonumber \\ 
	&+&(a^2-b^2)(a-b)^2+20ab^3\Bigg) \geq \frac{5b^2}{4ac^2}, \nonumber \\
	    \frac{dc}{dt}&=& \frac{c}{16 a^2 b^2 c^2}\Bigg(30a^4+3(a^2-b^2)^2+3(a^2-b^2)(a-b)^2+30b^4\Bigg) \nonumber \\ 
	    &\geq& \frac{15(a^4+b^4)}{8a^2b^2c}.
	\end{eqnarray}
	Thus $b$ and $c$ are non-decreasing functions and a is non-increasing. We then get 
	$a_{0} \geq a \geq b \geq b_{0}$.
	
	There is one conserved quantity for this geometry: for all times, one has
	\begin{equation}
	bc\frac{da}{dt}+ac\frac{db}{dt}+\frac{2ab}{3}\frac{dc}{dt}=0.
	\end{equation}
	The integrating factor for this equation is $c^{-1/3} $, so by direct integration, this equation gives $a b c^{2/3}=k$ where $k$ is some positive constant. Using this equation and also the identity $(a^4+b^4)=(a^2-b^2)^2+4a^2b^2$,  we find
	\begin{equation}
	    c^2(t) \geq 15t+c_{0}^2, 
	\end{equation}  
	and
	\begin{equation}
	\frac{1}{b^2} \frac{db}{dt} \leq \frac{5}{4ac^2}, \qquad \frac{1}{b^2} \geq \frac{-5}{2k} \int \frac{1}{c^{4/3}} dt, \qquad b^2 \leq \frac{2k}{-(15t+c_{0}^2)^{1/3}+k^{'}}.
	\end{equation}

Let us calculate the differences $a-c$ and $b-c$
\begin{eqnarray}
    \frac{d(a-c)}{dt}&=& \frac{a-b}{16 a^2 b^2 c^2}\biggl(-105a^4-45a^3b-43a^2b^2-63ab^3  \nonumber \\
	&-&c\Bigl(63(a-b)^3+192ab(a-b)\Bigr)\biggr).
  \end{eqnarray}
Therefore $a-c$ is a non-increasing function since $a_{0} \geq a \geq b \geq b_{0}$, and all the terms in the parenthesis are negative. 
\begin{eqnarray}
   \frac{d(b-c)}{dt}
     	&=&\frac{a-b}{16 a^2 b^2 c^2}\Bigg( 63a^3b+43a^2b^2+45ab^3+105b^4 \nonumber \\
	&-&c\Bigl(63(a-b)^3+192ab(a-b)\Bigr)\Bigg).
    \end{eqnarray}
We can not say much about the $b-c$ equation: it is inconclusive because there are negative and positive terms in the parenthesis. From two figures \ref{Fig.ISOMF1} and \ref{Fig.ISOMF2} showing the behavior of the geometry under the flow, we see that $a$ and $b$  initially approach each other and stay there, and when that happens, $c$ also reaches a constant value, generating a flat geometry.	
\begin{figure}[h]
    \centering
    \includegraphics[width=120mm,scale=0.5]{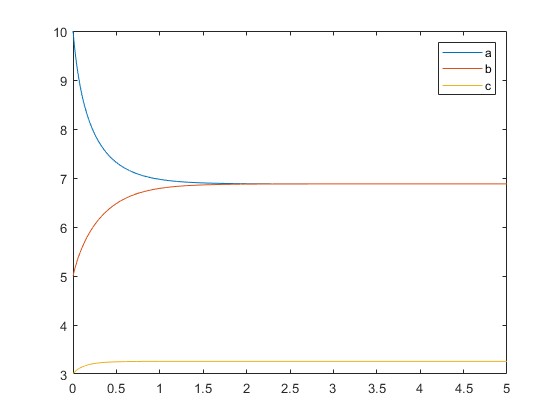}
    \caption{The behavior of  the $Isom(R^{2})$ geometry under the flow, first figure. The vertical axis represents the dimensions $a,b,c$. The horizontal axis represents the parameter of the flow $t$. Two dimensions $a,b$ approach each other and stay at a constant value, $c$ approaches to a constant value }
    \label{Fig.ISOMF1}
\end{figure}
\begin{figure}[h]
    \centering
    \includegraphics[width=120mm,scale=0.5]{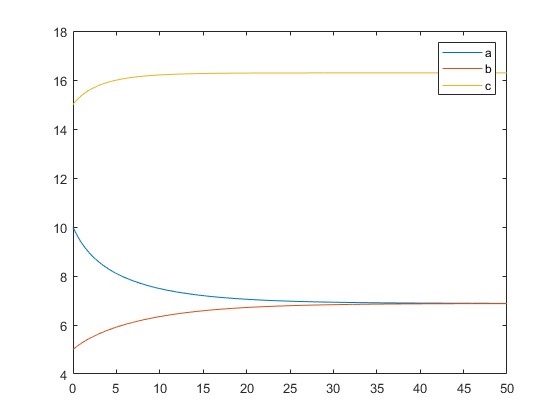}
    \caption{The behavior of  the $Isom(R^{2})$ geometry under the flow, second figure. The vertical axis represents the dimensions $a,b,c$. The horizontal axis represents a parameter of the flow $t$.}
    \label{Fig.ISOMF2}
\end{figure}

The curvature scalar when $a=b$ and $c$ is a constant, vanishes showing that the geometry is flat. At this point, $K$ also vanishes.
	\begin{equation}
		R=\frac{1}{2abc}(-a^2+2ab-b^2)=0,
	\end{equation}
\begin{equation}
K =\dfrac{1}{32a^2 b^2 c^2}(-21  a^4+20 a^3  b+2a^2  b^2+20a b^3-21 b^4)=0.
\end{equation}
\vspace {0.5 cm}

{\bf 5:  The geometry of $Sol$ with the structure constants $\lambda=-1 , \mu=0 , \nu=1 $}	
\vspace {0.5 cm}

For this geometry the flow equations are
		\begin{eqnarray} 
		\frac{da}{dt}&=&\frac{a}{16 a^2 b^2 c^2}\Big(-105a^4-60a^3c+2a^2c^2+20ac^3+63c^4\Big),  \nonumber \\
		\frac{db}{dt}&=& \frac{b}{16 a^2 b^2 c^2}\Big(63a^4+60a^3c-6a^2c^2+60ac^3+63c^4\Big), \nonumber \\
		\frac{dc}{dt}&=& \frac{c}{16 a^2 b^2 c^2}\Big(63a^4+20a^3c+2a^2c^2-60ac^3-105c^4\Big).
	\end{eqnarray}
	These equations are symmetric in $a$ and $c$ so we look at the difference between these two dimensions,
	\begin{equation}
	    	\frac{d(a-c)}{dt}=-\frac{3(a-c)}{16 a^2 b^2 c^2}\Big((a+c)^2(35a^2+6ac+35c^2)\Big).
	\end{equation}
		If at some time $t=\tau$, $a=c$ then the flow stops, so if initially $a_{0}\geq c_{0}$ then for all times $a \geq c$ and this implies also that $(a-c)$ is a non-increasing function. Using this condition on $a$ and $c$, we can make some estimates about the behavior of the flow equations. After some manipulations, we get
			\begin{eqnarray} 
		\frac{da}{dt}&=&\frac{a}{16 a^2 b^2 c^2}\Big(-80a^4-60c(a^3-c^3)-2a^2(a^2-c^2)  \nonumber \\
		&-&20a(a^3-c^3)-3(a^4-c^4)\Big) \leq -\frac{5a^3}{b^2c^2},  \nonumber  \\
		\frac{db}{dt}&=& \frac{b}{16 a^2 b^2 c^2}\Big(60a^4+60a^3c+3(a^2-c^2)^2+60ac^3+60c^4\Big)  \nonumber \\ 
		&\geq & \frac{15b}{4 a^2 b^2 c^2}\left(a^4+a^3c+ac^3+c^4\right),  \nonumber  \\
		\frac{dc}{dt}&=& \frac{c}{16 a^2 b^2c^2}\left(63a^4+20a^3c+2a^2c^2-60ac^3-105c^4\right).   
	\end{eqnarray}
	Thus $a$ is decreasing, and $b$ is increasing, but the behavior of $c$ is inconclusive from this construction. We write $a_{0} \geq a \geq c$. There is one conserved quantity, namely
	\begin{equation}
	    bc\frac{da}{dt}+ac\frac{2}{3}\frac{db}{dt}+ab\frac{dc}{dt}=0.
	\end{equation}
	The integrating factor for this equation is $b^{-1/3} $ , so we have $ab^{2/3}c=k=a_{0}b_{0}^{2/3}c_{0}$. 
	
	From the figure \ref{Fig.SOL}, we see that $a$ and $c$ approach each other and then go to a constant value, while $b$ always increases generating a cigar degeneracy.
\begin{figure}[h]
    \centering
    \includegraphics[width=110mm,scale=0.5]{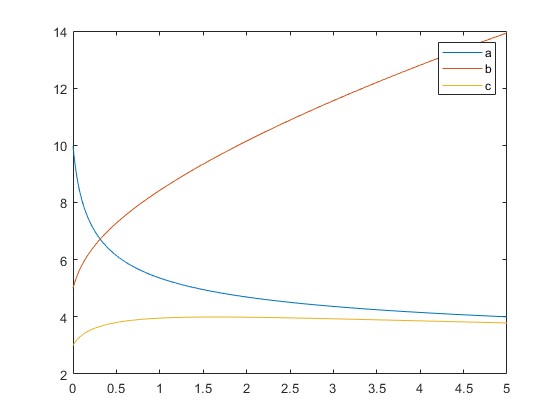}
    \caption{The behavior of  the $Sol$ geometry under the flow. The vertical axis represents the dimensions $a,b,c$. The horizontal axis represents the parameter of the flow $t$. Two dimensions $a,c$ first approach each other then shrink, $b$ expands. }
    \label{Fig.SOL}
\end{figure}

When $a=c$ the flow equations show that $a,c$ will eventually shrink, but are controlled by the constant $k$. 
		\begin{eqnarray} 
		\frac{da}{dt}&=&-\frac{5a}{b^2}, \qquad \frac{db}{dt}= \frac{15}{ b}, \nonumber \\
		 b^{2}(t)&=& b^{2}_{0}+30t, \qquad a(t)=a_{0}b_{0}^{1/3} b^{-1/3}(t).	
	\end{eqnarray}
The scalar curvature and the trace $K$ then go to zero:
\begin{equation}
		R=\frac{1}{2abc}(-a^2-2ac-c^2)=-\frac{2}{b}, \nonumber
	\end{equation}
	\begin{eqnarray}
K =\dfrac{1}{32a^2 b^2 c^2}\left(-21  a^4-20  a^3  c +2a^2  c^2-20 a  c^3-21 c^4\right)=-\dfrac{5}{2 b^2 }.
\end{eqnarray}
\vspace {0.5 cm}
{\bf 6: The geometry of $Nil$ with the structure constants $\lambda=-1 , \mu=0 , \nu=0 $}
\vspace {0.5 cm}

The flow equations for this geometry are
			\begin{eqnarray} 
		\frac{da}{dt}=-\frac{105a^3}{16 b^2 c^2}, \qquad 
		\frac{db}{dt}= \frac{63a^2b}{16 b^2 c^2}, \qquad
		\frac{dc}{dt}= \frac{63a^2c}{16 b^2 c^2}.
	\end{eqnarray}
	These equations are symmetric in $b$ and $c$. There is one conserved quantity for this geometry
	\begin{equation}
	    \frac{6}{5}\frac{1}{a}\frac{da}{dt}+\frac{1}{b}\frac{db}{dt}+\frac{1}{c}\frac{dc}{dt}=0,
	\end{equation}
	hence we have $a^{6/5}bc=k=a_{0}^{6/5}b_{0}c_{0}$. Furthermore, we have
	\begin{equation}
	    \frac{d}{dt}\left(\frac{b^2c^2}{a^2}\right)=\frac{231}{8}, \qquad \frac{b^2c^2}{a^2}=\frac{231}{8} t + \frac{b_{0}^2c_{0}^2}{a_{0}^2}.
	\end{equation}
	We find by integration
	\begin{eqnarray}
	 a(t)&=&a_{0}\left(\frac{b_{0}^2c_{0}^2}{\frac{231}{8}t+\frac{b_{0}^2c_{0}^2}{a_{0}^2}}\right)^{5/22}, \quad
	 b(t)=b_{0}\left(\frac{a_{0}^{12/5}c_{0}^2}{b_{0}^{14/5}}\right)\left(\frac{231}{8} t + \frac{b_{0}^2c_{0}^2}{a_{0}^2}\right)^{3/22}, \nonumber \\
	 c(t)&=&c_{0}\left(\frac{a_{0}^{12/5}b_{0}^2}{c_{0}^{14/5}}\right)\left(\frac{231}{8} t + \frac{c_{0}^2b_{0}^2}{a_{0}^2}\right)^{3/22}.
	\end{eqnarray}
Thus while one dimension shrinks to a point for large $t$ as  $\sim t^{-5/22}$, the other two dimensions expand with $\sim t^{3/22}$ generating a  pancake degeneracy. 
\begin{figure}[h]
    \centering
    \includegraphics[width=120mm,scale=0.5]{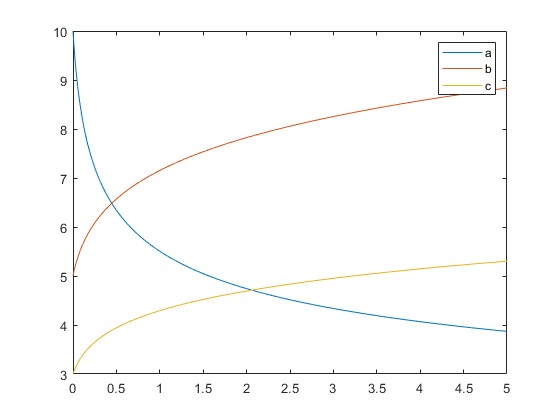}
    \caption{The behavior of  the $Nil$ geometry under the flow. The vertical axis represents the dimensions $a,b,c$. The horizontal axis represents the parameter of the flow $t$. Two dimensions $b,c$ expand, $c$ vanishes.}
    \label{Fig.NIL}
\end{figure}
	
The scalar curvature and trace $K$ go to zero as seen from the equations. 
\begin{eqnarray}
	R&=&-\frac{a}{2bc}=-\frac{1}{2}\Biggl(\frac{231}{8} t + \frac{b_{0}^2c_{0}^2}{a_{0}^2}\Biggr)^{-1/2}, \nonumber \\
	K &=&-\dfrac{21  a^2}{32 b^2 c^2}=-\frac{21}{32}\Biggl(\frac{231}{8} t + \frac{b_{0}^2c_{0}^2}{a_{0}^2}\Biggr)^{-1}.
\end{eqnarray}

\subsection{ The $K$ Tensor under the flow}

In this section, we examine the evolution of the basic geometric tensors and of the $K$ tensor itself under the $ K$-flow.  The following computations are similar to those of variations in an orthonormal frame background. We still study these equations in normal coordinates where at a point all partial derivatives of the metric $\partial_{k} g$ and Christoffel symbols $\Gamma$ vanish, but their time derivatives, $\partial_{k}\partial_{t} g$ and $\partial_{t}\Gamma$ do not, and we also use the covariant conservation of $g$ and $K$. We begin with the assumption that we have a time-dependent Riemannian  metric that changes in time according to the proposed flow
\begin{equation}
\partial_{t}g_{ij}(t)= \alpha K_{ij}(t), \qquad \partial_{t}g^{ij}(t)= -\alpha K^{ij}(t) .
\end{equation}
The metric determinant flows as
	\begin{equation}
		\partial_{t}g=\frac{\delta g}{\delta g_{ij}}\partial_{t}g_{ij}=gg^{ij}\alpha K_{ij}=\alpha gK,
	\end{equation}
while the volume element flows as 
	\begin{equation}
	\partial_{t}d\mu=\partial_{t}\sqrt{g}dx^{1}\wedge dx^2 \wedge dx^{3} =\frac{1}{2}\alpha Kd\mu.
	\end{equation}
The antisymmetric Levi-Civita tensor reads
	\begin{equation}
	\partial_{t}\eta^{ijk}=\partial{t}\frac{\epsilon^{ijk}}{\sqrt{g}}=-\frac{1}{2}\alpha K \eta^{ijk}.
	\end{equation}
The Christoffel symbols flow as
\begin{eqnarray}
\partial_{t}\Gamma^{i}_{jk}(t)=\alpha \frac{1}{2}g^{il}\left(-\nabla_{l}K_{jk}+\nabla_{j}K_{kl}+\nabla_{k}K_{jl}\right).
\end{eqnarray}
The Riemann tensor flows as 
	\begin{eqnarray}
	\partial_{t}R_{\phantom{i}jkl}^{i}(t)&=&\partial_{t}\partial_{k}\Gamma^{i}_{jl}+\partial_{t}\Gamma^{i}_{kn}\Gamma^{n}_{jl}+\Gamma^{i}_{kn}\partial_{t}\Gamma^{n}_{jl} -k\rightleftharpoons l \nonumber \\
&=&\partial_{k}\partial_{t}\Gamma^{i}_{jl}-\partial_{l}\partial_{t}\Gamma^{i}_{jk}=\nabla_{k}\partial_{t}\Gamma^{i}_{jl}-\nabla_{l}\partial_{t}\Gamma^{i}_{jk}\nonumber \\
&=&\frac{1}{2}\alpha \left(-\nabla_{k}\nabla^{i}K_{jl}+\nabla_{k}\nabla_{j}K^{\phantom{l}i}_{l}+\nabla_{k}\nabla_{l}K^{\phantom{j}i}_{j}\right) \nonumber  \\
&+&\frac{1}{2}\alpha\left(\nabla_{l}\nabla^{i}K_{jk}-\nabla_{l}\nabla_{j}K^{\phantom{k}i}_{k}-\nabla_{l}\nabla_{k}K^{\phantom{j}i}_{j}\right). 
	\end{eqnarray}
The Ricci tensor flows as
\begin{eqnarray}
\partial_{t}R_{ij}(t)&=&\partial_{t}R_{\phantom{k}ikj}^{k}(t)=\frac{1}{2}\alpha \left(-\Delta K_{ij}+\nabla_{k}\nabla_{i}K^{\phantom{j}k}_{j}+\nabla_{k}\nabla_{j}K^{\phantom{i}k}_{i}-\nabla_{i}\nabla_{j}K\right) \nonumber \\
&=&\alpha\left(-g_{ij}S_{kl}K^{kl}+\frac{3}{2}(S^{\phantom{i}k}_{i}K_{kj}+K^{\phantom{i}k}_{i}S_{kj})+SK_{ij}-KS_{ij}\right) \nonumber \\
&-&\frac{1}{2}\alpha(\Delta K_{ij}+\nabla_{i}\nabla_{j}K). 
\end{eqnarray}
In the first line, we commuted the covariant derivatives and used the covariant conservation of the $K$ tensor, $\nabla_{k}K^{\phantom{i}k}_{i}=0$. 
The curvature scalar flows as 
	\begin{equation}
		\partial_{t}R(t)=\partial_{t}(g^{ij}R_{ij})=-\alpha K^{ij}R_{ij}+g^{ij}\partial_{t}R_{ij}=
	-\alpha \left (S_{kl}K^{kl}+SK+\Delta K\right).
	\end{equation}
The Schouten tensor flows as
\begin{eqnarray}
\partial_{t}S_{ij}(t)&=&\alpha\left(\frac{1}{4}g_{ij}(-3S_{kl}K^{kl}+SK+\Delta K)+\frac{3}{2}(S^{\phantom{i}k}_{i}K_{kj}+K^{\phantom{i}k}_{i}S_{kj})-KS_{ij}\right) \nonumber \\
&-&\frac{1}{2}\alpha(\Delta K_{ij}+\nabla_{i}\nabla_{j}K).
\end{eqnarray}
For a generic symmetric (0,2) tensor $T$ (we set $\alpha=-1$ from now on) one has 
\begin{eqnarray}
\partial_{t}\nabla_{k}T_{ij}&=&\partial_{t}(\partial_{k}T_{ij}-\Gamma^{n}_{ki}T_{nj}-\Gamma^{n}_{kj}T_{in})=\partial_{k}\partial_{t}T_{ij}-\partial_{t}(\Gamma^{n}_{ki}T_{nj})-\partial_{t}(\Gamma^{n}_{kj}T_{in}) \nonumber \\
&=&\nabla_{k}\partial_{t}T_{ij}-\partial_{t}(\Gamma^{n}_{ki})T_{nj}-(\partial_{t}\Gamma^{n}_{kj})T_{in},
\end{eqnarray}
and 
\begin{eqnarray}
\partial_{t}\nabla_{l}\nabla_{k}T_{ij}&=&\partial_{t}(\partial_{l}\nabla_{k}T_{ij}-\Gamma^{n}_{lk}\nabla_{n}T_{ij}-\Gamma^{n}_{li}\nabla_{k}T_{nj}-\Gamma^{n}_{lj}\nabla_{k}T_{in}) \nonumber \\
&=&\nabla_{l}\partial_{t}\nabla_{k}T_{ij}-\partial_{t}(\Gamma^{n}_{lk})\nabla_{n}T_{ij}-\partial_{t}(\Gamma^{n}_{li})\nabla_{k}T_{nj}-\partial_{t}(\Gamma^{n}_{lj})\nabla_{k}T_{in} \nonumber \\
&=&\nabla_{l}\nabla_{k}\partial_{t}T_{ij}-(\nabla_{l}\partial_{t}\Gamma^{n}_{ki})T_{nj}-(\nabla_{l}\partial_{t}\Gamma^{n}_{kj})T_{in}-\partial_{t}\Gamma^{n}_{ki}\nabla_{l}T_{nj}
 \nonumber \\
&-&\partial_{t}\Gamma^{n}_{kj}\nabla_{l}T_{in}-\partial_{t}\Gamma^{n}_{lk}\nabla_{n}T_{ij} 
-\partial_{t}\Gamma^{n}_{li}\nabla_{k}T_{nj}-\partial_{t}\Gamma^{n}_{lj}\nabla_{k}T_{in}.
\end{eqnarray}
The 3D ``Bach''-type tensor as defined in (\ref{Bach}) flows as
\begin{align}
\partial_{t}H^{ij}=&  \frac{K}{2}H^{ij} - \frac{1}{2}K_{lm}\eta^{ikl}\nabla_{k}C^{mj}-\frac{1}{2}\left(  \frac{1}{4}\nabla^{i}\nabla^{j}(-[SK]+SK+\Delta K)+\frac{1}{2}\Delta^{2}K^{ij} \right. \nonumber \\
&\left. +\Delta (KS^{ij}-2S^{ik}K^{\phantom{k}j}_{k})+ \frac{1}{4}g^{ij}\Delta (3[SK]-SK-\Delta K) \right.  \nonumber\\
&\left. + g^{ij}(\frac{1}{2}S^{kl}\Delta K_{kl}-\frac{1}{2}S^{kl}\nabla_{k}\nabla_{l} K -\frac{1}{2}\nabla_{k}S\nabla^{k}K+K[SS]-2[SKS])\right.  \nonumber \\
&\left.  -\frac{5}{2}S^{ik}\Delta K^{\phantom{k}j}_{k}  -\frac{1}{2}S\Delta K^{ij}-3 S^{ij}[SK]+K^{ij}[SS]+  S^{ik} K^{lj} S_{kl} +5 S^{ik} S^{lj} K_{kl} \right.  \nonumber \\
&\left. -4 S^{ik}  S^{\phantom{k}j}_{k} K+3S^{ik}  K^{\phantom{k}j}_{k} S +\frac{1}{2}S^{kl}\nabla_{k}\nabla_{l}K^{ij}-2S_{kl}\nabla^{k}\nabla^{i}K^{lj}  +\frac{1}{2}  S_{kl}\nabla^{i}\nabla^{j}K^{kl} \right.   \nonumber\\
&\left.+   S^{kj}\nabla_{k}\nabla^{i}K+2\nabla^{i}K_{kl}\nabla^{k}S^{lj} - 2\nabla^{k}K^{il}\nabla_{l}S^{\phantom{k}j}_{k}-\nabla_{k}S^{ij}\nabla^{k}K   +   \nabla^{i}S^{kj}\nabla_{k}K \right.  \nonumber \\                                               
&\left. +\frac{1}{2}\nabla^{k}K^{ij}\nabla_{k}S-\nabla^{i}K^{kj}\nabla_{k}S+K_{kl} \nabla^{k}\nabla^{i}S^{lj}-K^{kj} \nabla_{k}\nabla^{i}S -\frac{1}{2}K\nabla^{i}\nabla^{j}S \right) \nonumber \\
&+ i \leftrightarrow j,
\end{align}
where we set $[SK]:=S_{ij}K^{ij}, \quad [SKS]:=S^{\phantom{i}j}_{i}K^{\phantom{j}k}_{j}S^{\phantom{k}l}_{k}$ etc. for simplicity.  This tedious computation can be checked in various ways, one of them is to observe that, since $H=0$,  $\partial_{t}H=0$ follows from the above equation as expected. 
Similarly, the $J$-tensor as defined in (\ref{Jtensor}) flows as 
\begin{align}
\partial_{t}J^{ij}=&  \frac{3K}{2}J^{ij} -\frac{1}{2}\left( S ^{ij}(-\frac{9}{4}[SK]-\frac{1}{4}SK+\frac{3}{4}\Delta K) -S^{ik}(\Delta K^{\phantom{k}j}_{k}+\nabla_{k} \nabla^{j}K ) \right.  \nonumber\\
&\left. +g^{ij}\Bigl(\frac{9}{4}S[SK]+\frac{1}{4}S^{2}K-\frac{3}{4}S\Delta K -3[SKS])+ \frac{1}{2}S^{kl}(\Delta K_{kl}+\nabla_{k}\nabla_{l}K)\Bigr) \right.  \nonumber\\
&\left.  + \frac{1}{2}S(\Delta K^{ij}+\nabla^{i}\nabla^{j} K) -3S^{ik}  K^{\phantom{k}j}_{k} S +3S^{ik}S^{lj}K_{kl}+3S^{ik}K^{lj}S_{kl}  \right)   \nonumber \\
&+ i \leftrightarrow j,
\end{align}
while its trace flows as 
\begin{align}
\partial_{t}J&=  -K _{ij}J ^{ij}+3KJ+3[SKS]- \frac{3}{2}S[SK]-\frac{1}{2}S^{2}K+\frac{1}{2}S\Delta K  \nonumber\\
&-\frac{1}{2}S ^{ij}(\Delta K _{ij}+\nabla_{i}\nabla_{j}K).
\end{align}
Finally, we can write the flow of the  $K$-tensor:
\begin{eqnarray}
\partial_{t}K^{ij}&=&  \frac{3}{2}KK^{ij}+2KJ^{ij}- K_{lm}\eta^{ikl}\nabla_{k}C^{mj}+\frac{S^{ij}}{4}(21[SK]+SK-3\Delta K) \nonumber \\
&+&g^{ij}\Bigl( 5[SKS]-\frac{9}{4}S[SK]-K[SS]-\frac{1}{4}S^{2}K-\frac{3}{4}\Delta ([SK])+\frac{1}{4}\Delta (SK) \nonumber \\
&+&\frac{3}{4}S\Delta K +\frac{1}{4}\Delta^{2} K-S^{kl}\Delta K_{kl}+\frac{1}{2}\nabla^{k}S\nabla_{k}K \Bigr) -S^{2}K^{ij}-8S^{ik}S^{lj}K_{kl}\nonumber \\
&-&4S^{ik}K^{lj}S_{kl}+4S^{ik}  S^{\phantom{k}j}_{k}K+\frac{7}{2}S^{ik}\Delta K^{\phantom{k}j}_{k}-\frac{1}{2}S\nabla^{i}\nabla^{j}K-\frac{1}{2}S^{kl}\nabla_{k}\nabla_{l}K^{ij} \nonumber  \\
&+&2 S_{kl}\nabla^{k}\nabla^{i}K^{lj}- \frac{1}{2}S^{kl}\nabla^{i}\nabla^{j}K_{kl}-2\nabla^{i}K_{kl}\nabla^{k}S^{lj}+2\nabla^{k}K^{il}\nabla_{l}S^{\phantom{k}j}_{k}+\nabla_{k}S^{ij}\nabla^{k}K   \nonumber \\                                               
&-& \nabla^{i}S^{kj}\nabla_{k}K-\frac{1}{2}\nabla^{k}K^{ij}\nabla_{k}S+\nabla^{i}K^{kj}\nabla_{k}S-K_{kl} \nabla^{k}\nabla^{i}S^{lj}+K^{kj} \nabla_{k}\nabla^{i}S   \nonumber \\
&+&\frac{1}{2}K\nabla^{i}\nabla^{j}S-\frac{1}{2}\nabla^{2}K^{ij}+\frac{1}{4}\nabla^{i}\nabla^{j}([SK]-SK-\Delta K)+\Delta(2S^{ik}K^{\phantom{k}j}_{k}-KS^{ij}) \nonumber \\
&+& i \leftrightarrow j,
\end{eqnarray}
and its trace
\begin{eqnarray}
\partial_{t}K&=&  -2K _{ij}J ^{ij}+3K^{2}+6[SKS]- 3S[SK]-S^{2}K+S\Delta K \nonumber \\
&-&S ^{ij}(\Delta K _{ij}+\nabla_{i}\nabla_{j}K).
\end{eqnarray}

\subsection{ Conformal Properties of the $K$ Tensor} 
	
	Let us study the effects of a conformal mapping defined as $\tilde{g}_{ij}=e^{2\phi}g_{ij}$  on the $K$-tensor. To carry out this rather long computation, we make a few notes. The Levi-Civita symbol (in 3 dimensions) $\epsilon^{ijk}$ is a tensor density of weight +1. So we define $\eta^{ijk}=\epsilon^{ijk}/\sqrt{g }$ as a true tensor since the metric $g$ has weight +2. Thus the Cotton tensor defined as $C^{ij}=\eta^{ikl}\nabla_{k}S^{\phantom{l}j}_{l}$ is not conformally invariant, but under conformal transformations, it transforms as
$\widetilde{C^{ij}}=e^{-5\phi}C^{ij}$. It is possible to define a conformally invariant Cotton tensor of weight +5/3 as $C^{'ij}=g^{5/6}\eta^{ikl}\nabla_{k}S^{\phantom{l}j}_{l}$, but we will continue with the non-conformal tensor definition. Then we write 
\begin{eqnarray}
&&\widetilde{C^{\phantom{j}i}_{j}}=e^{-3\phi}C^{\phantom{j}i}_{j},\qquad \widetilde{C_{ij}}=e^{-\phi}C_{ij},\nonumber \\
&&\widetilde{\eta^{ijk}}=e^{-3\phi}\eta^{ijk},\quad \sqrt{\tilde{g} }=e^{3\phi}\sqrt{g },\quad  \widetilde{\eta^{\phantom{i}jk}_{i}}=e^{-\phi}\eta^{\phantom{i}jk}_{i} ,\nonumber \\
&&\widetilde{\nabla_{k}C_{ij}}=e^{-\phi}(\nabla_{k}C_{ij}-2C_{ij}\nabla_{k}\phi).
\end{eqnarray}
Under the conformal transformations, the 3D Bach tensor (\ref{Bach}) transforms as  
\begin{equation}
 \widetilde{H_{ij}}=e^{-2\phi}H_{ij}-2e^{-2\phi}\eta^{\phantom{i}kl }_{i}C_{lj}\nabla_{k}\phi,
\end{equation}
while the conformal transformation of the $J$ tensor (\ref{Jtensor}) is
\begin{eqnarray}
 \widetilde{J_{ij}}&=&e^{-2\phi}J_{ij}-\frac{e^{-2\phi} }{2} \Biggl( g_{ij} \Bigl( 2S_{kl}\nabla^{k} \nabla^{l} \phi -2S_{kl}\nabla^{k}\phi \nabla^{l} \phi -2S\Delta \phi +S\nabla_{k}\phi \nabla^{k}\phi \nonumber \\
&+& \Delta \phi\Delta \phi-\Delta \phi \nabla_{k}\phi \nabla^{k}\phi-\nabla_{k} \nabla_{l} \phi\nabla^{k} \nabla^{l} \phi+2 \nabla_{k} \nabla_{l} \phi\nabla^{k}\phi \nabla^{l} \phi \nonumber \\
&-&\frac{1}{2}\nabla_{k}\phi \nabla^{k}\phi\nabla_{l}\phi \nabla^{l}\phi \Bigr)+S_{ij}\Bigl( 2\Delta \phi    -	\nabla_{k}\phi \nabla^{k}\phi	\Bigr) +S\Bigl(2\nabla_{i} \nabla_{j}\phi -2\nabla_{i}\phi \nabla_{j}\phi	\Bigr) \nonumber \\
&+&S_{ik}\Bigl(-2\nabla_{j} \nabla^{k}\phi +2\nabla_{j}\phi \nabla^{k}\phi	\Bigr)+ S_{jk}\Bigl(-2\nabla_{i} \nabla^{k}\phi +2\nabla_{i}\phi \nabla^{k}\phi	\Bigr) \nonumber \\
&+&	\nabla_{i} \nabla_{j}\phi \Bigl( -2   \Delta \phi + \nabla_{k}\phi \nabla^{k}\phi    \Bigr) 	+\nabla_{i}\phi  \nabla_{j}\phi \Bigl( 2   \Delta \phi + \nabla_{k}\phi \nabla^{k}\phi    \Bigr) \nonumber \\
&+&2\nabla_{i}\nabla_{k}\phi \nabla^{k}\nabla_{j}\phi -2\nabla_{i}\nabla_{k}\phi \nabla^{k}\phi\nabla_{j}\phi-2\nabla_{j}\nabla_{k}\phi \nabla^{k}\phi \nabla_{i}\phi  \Biggr).
\end{eqnarray}
From the last two equations, one can find the conformal transformations of the $K_{ij}$ tensor.

\section{Acknowledgments}
We would like to thank J.D. Streets for pointing out his work and solutions on the existence-uniqueness issues, and Sahin Ulas Koprucu for his help with the \textsc{Matlab} code.

\end{document}